\documentclass[12pt]{article}
\usepackage{amsmath,amsfonts,latexsym,amsthm,amssymb}
\newcommand{\D}{\mathbb D_t^{(\alpha )}}
\newcommand{\Rn}{\mathbb R^n}
\numberwithin{equation}{section}

\begin{document}
\newtheorem{prop}{Proposition}
\newtheorem{teo}{Theorem}
\pagestyle{plain}
\title{Cauchy Problem for Fractional Diffusion-Wave Equations with Variable Coefficients}
\author{Anatoly N. Kochubei
\\ \footnotesize Institute of Mathematics,\\
\footnotesize National Academy of Sciences of Ukraine,\\
\footnotesize Tereshchenkivska 3, Kiev, 01601 Ukraine\\
\footnotesize E-mail: \ kochubei@i.com.ua}
\date{}
\maketitle

\vspace*{3cm}
\begin{abstract}
We consider an evolution equation with the Caputo-Dzhrbashyan fractional derivative of order $\alpha \in (1,2)$ with respect to the time variable, and the second order uniformly elliptic operator with variable coefficients acting in spatial variables. This equation describes the propagation of stress pulses in a viscoelastic medium. Its properties are intermediate between those of parabolic and hyperbolic equations. In this paper,
we construct and investigate a fundamental solution of the Cauchy problem, prove existence and uniqueness theorems for such equations.
\end{abstract}
\vspace{2cm}


\newpage
\section{Introduction}

The fractional diffusion-wave equation has the form
\begin{equation}
\left( \D u\right) (t,x)-\Delta u(t,x)=f(t,x),\quad t\in (0,T],\ x\in
\mathbb R^n,
\end{equation}
where $1<\alpha <2$, $\D$ is the Caputo-Dzhrbashyan fractional derivative, that is
$$
\left( \mathbb D_t^{(\alpha )}u\right) (t,x)=\frac{1}{\Gamma (2-\alpha )}\frac{\partial}{\partial t}\int\limits_0^t (t-\tau )^{-\alpha +1}u'_\tau (\tau ,x)\,d\tau -t^{-\alpha +1}\frac{u'_t(0,x)}{\Gamma (2-\alpha )}.
$$

This equation describes the propagation of stress pulses in a viscoelastic medium \cite{M}; its properties are intermediate between those of the classical heat and wave equations. On the one hand, regularity properties of its solutions resemble those of parabolic equations; that follows, for example, from the representation of solutions as convolutions with Green kernels which are ordinary functions (possessing, if $n>1$, singularities with respect to the spatial variables, just as the fractional diffusion equations of order $\alpha \in (0,1)$ \cite{K90,EK}).

On the other hand, the well-posed Cauchy problem for the equation (1.1) requires two initial functions
\begin{equation}
u(0,x)=u_0(x),\quad u'_t(0,x)=u_1(x),
\end{equation}
as for the wave equation. The fundamental solution of the Cauchy problem (1.1)-(1.2) decays exponentially outside the ``fractional light cone'' $\{ |x|t^{-\alpha /2}\le 1\}$, which is the characteristic property of the class of fractional-hyperbolic equations and systems \cite{K13}.

The fractional diffusion-wave equation and its generalizations (linear and nonlinear fractional operator-differential equations with the operator $\mathbb D^{(\alpha )},1<\alpha <2$) have been studied by many authors; see \cite{CGL,CLS,Fu,Ha,LMP,LP,M96,MP,PP,Ps,SW,VZ} and references therein. Such equations can be interpreted also as special cases of abstract Volterra equations \cite{Pr}.

In this paper we consider equations of the form
\begin{equation}
Lu\equiv \left( \D u\right) (t,x)-\mathcal B u(t,x)=f(t,x)
\end{equation}
where $1<\alpha <2$,
$$
\mathcal Bu(t,x)=\sum\limits_{i,j=1}^na_{ij}(x)\frac{\partial^2 u(t,x)}{\partial
x_i\partial x_j}+\sum\limits_{j=1}^nb_j(x)\frac{\partial u(t,x)
}{\partial x_j}+c(x)u(t,x),
$$
and there exists such a constant $\delta_0>0$ that for any $x,\xi \in \Rn$
\begin{equation}
\sum\limits_{i,j=1}^na_{ij}(x)\xi_i\xi_j\ge \delta_0|\xi |^2.
\end{equation}
We assume that $a_{ij}=a_{ji}$, the coefficients $a_{ij},b_k,c$ are bounded, uniformly H\"older continuous real-valued functions with the H\"older exponent $\gamma$ satisfying the inequality
\begin{equation}
2-\frac2\alpha <\gamma \le 1.
\end{equation}
Conditions on the right-hand side $f$ and the initial functions $u_0,u_1$ will be stated below.

The main task is to construct and study the fundamental solution of the Cauchy problem for the equation (1.3), thus extending the classical Levi method well known for parabolic equations and systems (see, for example, \cite{E,EIK,F}). For fractional diffusion equations (the case $0<\alpha <1$), this method was implemented by Eidelman and the author \cite{EK}; see also \cite{EIK,K12}. Shortly before his death, S. D. Eidelman (1920--2005), one of the founders of general theory of parabolic equations and systems, discussed with the author a possibility to consider the case $1<\alpha <2$. At that time it looked quite difficult because the Levi method requires precise estimates of various kernels based on the fundamental solution of the equation with constant coefficients. A representation of this fundamental solution found in \cite{SW} and \cite{K90} and used in \cite{EK} involves Fox's H-function. In many cases, we had to use several terms of the asymptotic expansions of different H-functions appearing in a complicated expression, and to check that some terms are cancelled, in order to obtain the required estimates. A representation in terms of a more accessible Wright function was known only for $n=1$ \cite{M96,MP}.

This situation was changed in 2009 by the paper \cite{Ps} by Pskhu who found and investigated in detail an expression of a fundamental solution of the problem (1.1)-(1.2) (and its analogs for other notions of a fractional derivative) in terms of the Wright function. Our investigation of the equation (1.3) is based on the results from \cite{Ps}. Here and below the fundamental solution of the Cauchy problem is understood as a collection of three kernels $Z_1(t,x; \xi )$, $Z_2(t,x; \xi )$, $Y(t,x; \xi )$, such that the function
\begin{equation}
u(t,x)=\int\limits_{\mathbb R^n} Z_1(t,x;\xi )u_0(\xi )\,d\xi +\int\limits_{\mathbb R^n} Z_2(t,x;\xi )u_1(\xi )\,d\xi +\int\limits_0^td\tau \int\limits_{\mathbb R^n} Y(t-\tau,x,;\xi)f(\tau ,\xi)\,d\xi
\end{equation}
is, under some conditions upon $u_0,u_1,f$, a classical solution of the Cauchy problem. This means that

\begin{description}
\item[(i)] $u(t,x)$ is twice continuously differentiable in $x$
for each $t>0$;
\item[(ii)] for each $x\in \mathbb R^n$ $u(t,x)$ is continuously differentiable in
$(t,x)$ on $[0,T]\times \Rn$, and the fractional integral
$$
\left( I_{0+}^{2-\alpha }u\right) (t,x)=\frac{1}{\Gamma (2-\alpha
)}\int\limits_0^t(t-\tau )^{-\alpha +1}u'_\tau (\tau ,x)\,d\tau
$$
is continuously differentiable in $t$ for $t>0$;
\item[(iii)]
$u(t,x)$ satisfies the equation and initial conditions.
\end{description}

If the kernels in (1.6) depend on the difference $x-\xi$, and also on some parameter $\eta$, we will write them also as $Z_1(t,x-\xi ;\eta )$ etc. The iteration processes of the Levi method are carried out three times, separately for each of the above kernels.

The structure of this paper is as follows. In Section 2, we collect the results for equations with constant coefficients covered essentially by \cite{Ps}. In Section 3, we give a description of Levi's method for the case of diffusion-wave equations. Then, in Section 4 we consider the Cauchy problem (1.2)-(1.3), with a stress on some features different from both the classical theory of parabolic equations and the case of fractional diffusion equations. Section 5 is devoted to the uniqueness of a solution of the Cauchy problem.

Note that a complete exposition of the above material would be quite lengthy. Therefore we omit fragments of proofs identical to those appearing in the classical theory or in \cite{EK}. In such cases we give only formulations and the appropriate references. At the same time, the check of initial conditions and the proof of uniqueness require new assumptions and techniques, as compared to the case $0<\alpha <1$, and will be expounded in detail.

\section{Equations with Constant Coefficients}

{\bf 2.1. Constructions.} Let us begin with the case where the coefficients of $\mathcal B$ are constant, and only the leading terms are present, so that
$$
\mathcal B=\sum\limits_{i,j=1}^na_{ij}\frac{\partial^2}{\partial x_i\partial x_j}
$$
where $A=(a_{ij})$ is a positive definite symmetric matrix. In fact, we will need a little more general situation, in which the coefficients $a_{ij}$ depend on a parameter $\eta \in \Rn$, so that the functions $\eta \mapsto a_{ij}(\eta )$ are bounded and uniformly H\"older continuous, and the ellipticity condition (1.4) holds uniformly with respect to $\eta$.

Let $A(\eta )=(A^{(ij)})$ be the matrix inverse to $(a_{ij})$. Denote
$$
\mathfrak A(y,\eta )=\sum\limits_{i,j=1}^nA^{(ij)}(\eta )y_iy_j,\quad y,\eta \in \Rn.
$$
By our assumptions,
\begin{gather}
C_1|y|^2\le \mathfrak A(y,\eta )\le C_2|y|^2,\\
\left| \mathfrak A(y,\eta')-\mathfrak A(y,\eta'')\right| \le
C|\eta'-\eta''|^\gamma |y|^2,\\
\left| \left[ \det A(\eta')\right]^{1/2}-\left[ \det
A(\eta'')\right]^{1/2}\right| \le C|\eta'-\eta''|^\gamma.
\end{gather}
Here and below we denote by $C$ (with indices or exponents or without them) various positive constants. Positive constants appearing under the sign of exponential will be denoted $\sigma$, while the exponents of the H\"older continuity are all denoted by the same letter $\gamma$.

Given a fundamental solution $Z_{1,0}(t,x-\xi )$, $Z_{2,0}(t,x-\xi )$, $Y_0(t,x-\xi )$ of the Cauchy problem for the equation (1.1), we can write a similar triple for our present case setting
\begin{equation}
Z_k^{(0)}(t,x-\xi ;\eta )=\frac1{[\det A(\eta )]^{1/2}}Z_{k,0}(t,[\mathfrak A(x-\xi ,\eta )]^{1/2}),\quad k=1,2,
\end{equation}
\begin{equation}
Y^{(0)}(t,x-\xi ;\eta )=\frac1{[\det A(\eta )]^{1/2}}Y_0(t,[\mathfrak A(x-\xi ,\eta )]^{1/2});
\end{equation}
the motivation for this notation will become clear later. Note that the transition from $Z_{k,0}$ and $Y_0$ to $Z_k^{(0)}$ and $Y^{(0)}$ can be interpreted as a linear change of variables; see the proof of Theorem 1, Chapter 1, in \cite{F}.

A fundamental solution $(Z_{1,0},Z_{2,0},Y_0)$ was found by Pskhu \cite{Ps}:
\begin{equation}
Z_{1,0}(t,x)=D_{0t}^{\alpha -1}\Gamma_{\alpha ,n}(t,x),\quad Z_{2,0}(t,x)=D_{0t}^{\alpha -2}\Gamma_{\alpha ,n}(t,x),\quad Y_0(t,x)=\Gamma_{\alpha ,n}(t,x),
\end{equation}
where, following \cite{Ps}, we use the unified notation for the Riemann-Liouville integrals and derivatives: if $s$ is the initial point, and $p\in \mathbb N$, $p-1<\beta \le p$, then
$$
D_{st}^\beta g(t)=\operatorname{sign}^p (t-s)\left( \frac{\partial}{\partial t}\right)^pD_{st}^{\beta -p}g(t),
$$
$$
D_{st}^\mu g(t)=\frac{\operatorname{sign}(t-s)}{\Gamma (-\mu )}\int\limits_s^tg(\tau )(t-\tau )^{-\mu -1}d\tau ,\quad \mu <0,
$$
$D_{st}^0 g(t)=g(t)$. The function $\Gamma_{\alpha ,n}$ is defined as follows:
\begin{equation}
\Gamma_{\alpha ,n}(t,x)=c_n t^{\alpha -\frac{\alpha n}2-1}f_{\alpha /2}(t^{-\alpha /2}|x|;n-1,\alpha -\frac{\alpha n}2).
\end{equation}
Here
$$
f_{\alpha /2}(z;\mu ,\delta )=\begin{cases}
\frac2{\Gamma (\mu /2)}\int\limits_1^\infty \Phi (-\alpha /2,\delta ,-zt)(t^2-1)^{\frac{\mu}2-1}dt, &\text{if $\mu >0$};\\
\Phi (-\alpha /2,\delta ,-z), &\text{if $\mu =0$},\end{cases}
$$
$$
\Phi (-\alpha /2,\delta ,z)=\sum\limits_{m=0}^\infty \frac{z^m}{m!\Gamma (\delta -\frac{\alpha m}2)}
$$
is the Wright function, $c_n=2^{-n}\pi^{(1-n)/2}$. A series representation of the function $f_{\alpha /2}$ (see \cite{Ps}) shows that, if $n>1$, the fundamental solution has a singularity in spatial variables.

The work with expressions containing $f_{\alpha /2}$ is simplified by the identities \cite{Ps}
\begin{equation}
\frac{d}{dz}f_{\alpha /2}(z;\mu ,\delta )=-\frac{z}2 f_{\alpha /2}(z;\mu +2,\delta -\alpha );
\end{equation}
\begin{equation}
D_{st}^\zeta |t-s|^{\delta -1}f_{\alpha /2}(|t-s|^{-\alpha /2}z;\mu ,\delta )=|t-s|^{\delta -\zeta -1}f_{\alpha /2}(|t-s|^{-\alpha /2}z;\mu ,\delta -\zeta ), \quad \zeta \in \mathbb R.
\end{equation}

\medskip
{\bf 2.2. Estimates.} Using the estimates for the integer and fractional order derivatives of the function $\Gamma_{\alpha ,n}$ found in \cite{Ps} and the property (2.1), we find estimates of the kernels $Z_1^{(0)}, Z_2^{(0)}, Y^{(0)}$ and their derivatives (some higher derivatives absent in \cite{Ps} are treated easily using (2.8)). The estimates are different for $n\ge 3$, $n=2$, and $n=1$. Therefore we consider these cases separately. Denote
$$
\rho_\sigma (t,x,\xi )=\exp \left\{ -\sigma (t^{-\alpha /2}|x-\xi |)^{\frac2{2-\alpha}}\right\},\quad \sigma >0.
$$

Let $n\ge 3$. Then
\begin{equation}
\left| D_x^m Z_1^{(0)}(t,x-\xi ;\eta )\right| \le Ct^{-\alpha }|x-\xi |^{-n+2-|m|}\rho_\sigma (t,x,\xi ),\quad |m|\le 3;
\end{equation}
\begin{equation}
\left| D_x^m Z_2^{(0)}(t,x-\xi ;\eta )\right| \le Ct^{-\alpha +1}|x-\xi |^{-n+2-|m|}\rho_\sigma (t,x,\xi ),\quad |m|\le 3;
\end{equation}
\begin{equation}
\left| Y^{(0)}(t,x-\xi ;\eta )\right| \le Ct^{\alpha -\frac{\alpha n}2-1}\mu_n(t^{-\alpha /2}|x-\xi |)\rho_\sigma (t,x,\xi ),
\end{equation}
where
$$
\mu_n(z)=\begin{cases}
1, &\text{if $n=3$};\\
1+|\log z|, &\text{if $n=4$};\\
z^{-n+4}, &\text{if $n\ge 5$}.\end{cases}
$$

Next,
\begin{equation}
\left| \frac{\partial}{\partial x_i}Y^{(0)}(t,x-\xi ;\eta )\right| \le Ct^{-\alpha -1}|x-\xi |^{-n+3}\rho_\sigma (t,x,\xi ).
\end{equation}
Making the estimates a little rougher, we can unify (2.12) and (2.13), together with the estimates for second and third order derivatives, into the following unified estimate:
\begin{equation}
\left| D_x^m Y^{(0)}(t,x-\xi ;\eta )\right| \le Ct^{-1}|x-\xi |^{-n+2-|m|}\rho_\sigma (t,x,\xi ),\quad |m|\le 3,
\end{equation}
which will be used in the implementation of Levi's method. However the initial estimate (2.13) will also be useful (for the proof of the uniqueness theorem).

The above transformation of estimates is based on a procedure frequently used throughout the paper -- we can drop a positive power of the expression $t^{-\alpha /2}|x-\xi |$, simultaneously taking a smaller $\sigma >0$ in the factor $\rho_\sigma$.

The estimates for time derivatives of the functions $Z_1^{(0)}, Z_2^{(0)}, Y^{(0)}$ are as follows:
\begin{equation}
\left| \frac{\partial}{\partial t}Z_1^{(0)}(t,x-\xi ;\eta )\right| \le Ct^{-\alpha -1}|x-\xi |^{-n+2}\rho_\sigma (t,x,\xi );
\end{equation}
\begin{equation}
\left| \frac{\partial}{\partial t}Z_2^{(0)}(t,x-\xi ;\eta )\right| \le Ct^{-\alpha }|x-\xi |^{-n+2}\rho_\sigma (t,x,\xi );
\end{equation}
\begin{equation}
\left| \frac{\partial}{\partial t}Y^{(0)}(t,x-\xi ;\eta )\right| \le Ct^{\alpha -\frac{\alpha n}2-2}\mu_n(t^{-\alpha /2}|x-\xi |)\rho_\sigma (t,x,\xi ),
\end{equation}
\begin{equation}
\left| \D Z_1^{(0)}(t,x-\xi ;\eta )\right| \le Ct^{-2\alpha }|x-\xi |^{-n+2}\rho_\sigma (t,x,\xi );
\end{equation}
\begin{equation}
\left| \D Z_2^{(0)}(t,x-\xi ;\eta )\right| \le Ct^{-2\alpha +1}|x-\xi |^{-n+2}\rho_\sigma (t,x,\xi );
\end{equation}
\begin{equation}
\left| \D Y^{(0)}(t,x-\xi ;\eta )\right| \le Ct^{-\alpha -1}|x-\xi |^{-n+2}\rho_\sigma (t,x,\xi ).
\end{equation}

Let $n=2$. Then
\begin{equation}
\left| Z_1^{(0)}(t,x-\xi ;\eta )\right| \le Ct^{-\alpha }\left[ \left|\log \left( t^{-\alpha /2}|x-\xi |\right)\right| +1\right]\rho_\sigma (t,x,\xi );
\end{equation}
\begin{equation}
\left| Z_2^{(0)}(t,x-\xi ;\eta )\right| \le Ct^{-\alpha +1}\left[ \left|\log \left( t^{-\alpha /2}|x-\xi |\right)\right| +1\right]\rho_\sigma (t,x,\xi );
\end{equation}
\begin{equation}
\left| Y^{(0)}(t,x-\xi ;\eta )\right| \le Ct^{-1}\rho_\sigma (t,x,\xi );
\end{equation}
\begin{equation}
\left| D_x^mY^{(0)}(t,x-\xi ;\eta )\right| \le Ct^{-\alpha -1}\left[ \left|\log \left( t^{-\alpha /2}|x-\xi |\right)\right| +1\right]\rho_\sigma (t,x,\xi ),\quad |m|=1;
\end{equation}
\begin{equation}
\left| \frac{\partial}{\partial t}Z_1^{(0)}(t,x-\xi ;\eta )\right| \le Ct^{-\alpha -1}\left[ \left|\log \left( t^{-\alpha /2}|x-\xi |\right)\right| +1\right]\rho_\sigma (t,x,\xi );
\end{equation}
\begin{equation}
\left| \frac{\partial}{\partial t}Z_2^{(0)}(t,x-\xi ;\eta )\right| \le Ct^{-\alpha }\left[ \left|\log \left( t^{-\alpha /2}|x-\xi |\right)\right| +1\right]\rho_\sigma (t,x,\xi );
\end{equation}
\begin{equation}
\left| \frac{\partial}{\partial t}Y^{(0)}(t,x-\xi ;\eta )\right| \le Ct^{-2}\rho_\sigma (t,x,\xi );
\end{equation}
\begin{equation}
\left| \D Z_1^{(0)}(t,x-\xi ;\eta )\right| \le Ct^{-2\alpha }\left[ \left|\log \left( t^{-\alpha /2}|x-\xi |\right)\right| +1\right]\rho_\sigma (t,x,\xi );
\end{equation}
\begin{equation}
\left| \D Z_2^{(0)}(t,x-\xi ;\eta )\right| \le Ct^{-2\alpha +1}\left[ \left|\log \left( t^{-\alpha /2}|x-\xi |\right)\right| +1\right]\rho_\sigma (t,x,\xi );
\end{equation}
\begin{equation}
\left| \D Y^{(0)}(t,x-\xi ;\eta )\right| \le Ct^{-\alpha -1}\left[ \left|\log \left( t^{-\alpha /2}|x-\xi |\right)\right| +1\right]\rho_\sigma (t,x,\xi );
\end{equation}

The estimates (2.10) and (2.11) with $1\le m\le 3$, as well as the estimate (2.14) with $2\le m\le 3$, remain valid for $n=2$.

In the one-dimensional case ($n=1$), our kernels have no singularity with respect to the spatial variables. The estimates are as follows:
\begin{equation}
\left| D_x^m Z_1^{(0)}(t,x-\xi ;\eta )\right| \le Ct^{-\frac{\alpha}2(m+1)}\rho_\sigma (t,x,\xi );
\end{equation}
\begin{equation}
\left| D_x^m Z_2^{(0)}(t,x-\xi ;\eta )\right| \le Ct^{-\frac{\alpha}2(m+1)+1}\rho_\sigma (t,x,\xi );
\end{equation}
\begin{equation}
\left| D_x^m Y^{(0)}(t,x-\xi ;\eta )\right| \le Ct^{-\frac{\alpha}2(m-1)-1}\rho_\sigma (t,x,\xi ).
\end{equation}
Here $0\le m\le 3$. There is a more refined estimate for $m=1$:
\begin{equation}
\left| \frac{\partial}{\partial x}Y^{(0)}(t,x-\xi ;\eta )\right| \le Ct^{-1-\frac{\alpha}2}|x|\rho_\sigma (t,x,\xi ).
\end{equation}

Next,
\begin{equation}
\left| \frac{\partial}{\partial t}Z_1^{(0)}(t,x-\xi ;\eta )\right| \le Ct^{-1-\frac{\alpha}2}\rho_\sigma (t,x,\xi );
\end{equation}
\begin{equation}
\left| \frac{\partial}{\partial t}Z_2^{(0)}(t,x-\xi ;\eta )\right| \le Ct^{-\frac{\alpha}2}\rho_\sigma (t,x,\xi );
\end{equation}
\begin{equation}
\left| \frac{\partial}{\partial t}Y^{(0)}(t,x-\xi ;\eta )\right| \le Ct^{\frac{\alpha}2-2}\rho_\sigma (t,x,\xi );
\end{equation}
The estimates of $\D Z_1^{(0)}$, $\D Z_2^{(0)}$, $\D Y^{(0)}$ coincide with those for the appropriate second spatial derivatives.

\medskip
{\bf 2.3. Differences.} We will need estimates of values of $Z_1^{(0)}, Z_2^{(0)}$, and $Y^(0)$ at different values of the parameter $\eta$.

\medskip
\begin{prop}
For any $\eta',\eta''\in \Rn$, $t>0$, $x,\xi \in \Rn$, for each multi-index $m$, $|m|\le 2$,
\begin{equation}
\left| D_x^m Z_1^{(0)}(t,x-\xi ;\eta' )-D_x^m Z_1^{(0)}(t,x-\xi ;\eta'')\right| \le C|\eta'-\eta''|^\gamma t^{-\alpha }|x-\xi |^{-n+2-|m|}\rho_\sigma (t,x,\xi );
\end{equation}
\begin{equation}
\left| D_x^m Z_2^{(0)}(t,x-\xi ;\eta' )-D_x^m Z_2^{(0)}(t,x-\xi ;\eta'')\right| \le C|\eta'-\eta''|^\gamma t^{-\alpha +1}|x-\xi |^{-n+2-|m|}\rho_\sigma (t,x,\xi );
\end{equation}
\begin{equation}
\left| D_x^m Y^{(0)}(t,x-\xi ;\eta' )-D_x^m Y^{(0)}(t,x-\xi ;\eta'')\right| \le C|\eta'-\eta''|^\gamma t^{-1}|x-\xi |^{-n+2-|m|}\rho_\sigma (t,x,\xi );
\end{equation}
\begin{equation}
\left| \frac{\partial}{\partial t}Z_1^{(0)}(t,x-\xi ;\eta' )-\frac{\partial}{\partial t} Z_1^{(0)}(t,x-\xi ;\eta'')\right| \le C|\eta'-\eta''|^\gamma t^{-\alpha -1}|x-\xi |^{-n+2}\rho_\sigma (t,x,\xi );
\end{equation}
\begin{equation}
\left| \frac{\partial}{\partial t}Z_2^{(0)}(t,x-\xi ;\eta' )-\frac{\partial}{\partial t} Z_2^{(0)}(t,x-\xi ;\eta'')\right| \le C|\eta'-\eta''|^\gamma t^{-\alpha }|x-\xi |^{-n+2}\rho_\sigma (t,x,\xi );
\end{equation}
\begin{equation}
\left| \frac{\partial}{\partial t}Y^{(0)}(t,x-\xi ;\eta' )-\frac{\partial}{\partial t} Y^{(0)}(t,x-\xi ;\eta'')\right| \le C|\eta'-\eta''|^\gamma t^{-2}|x-\xi |^{-n+2}\rho_\sigma (t,x,\xi );
\end{equation}
\end{prop}

\medskip
{\it Proof}. Let us prove the inequality (2.38) with $m=0$. The proofs of all other estimates are similar.

Denote $\varphi_t(r)=D_{0t}^{\alpha -1}\Gamma_{\alpha ,n}(t,r)$ understanding $\Gamma_{\alpha ,n}(t,r)$ as the right-hand side of (2.7) with $|x|=r$. Using (2.7) we see that
$$
\frac{d}{dr}\varphi_t(r)=c_nD_{0t}^{\alpha -1}\left[ t^{\frac{\alpha }2-\frac{\alpha n}2-1}f'_{\alpha /2}(t^{-\alpha /2}r;n-1,\alpha -\frac{\alpha n}2)\right]
$$
where $f'_{\alpha /2}$ means the derivative with respect to the first argument. Then the identities (2.8) and (2.9) show that
\begin{multline*}
\frac{d}{dr}\varphi_t(r)=-\frac{c_nr}2 D_{0t}^{\alpha -1}\left[ t^{-\frac{\alpha n}2-1}f_{\alpha /2}(t^{-\alpha /2}r;n+1,-\frac{\alpha n}2)\right] \\
=-\frac{c_nr}2t^{-\frac{\alpha n}2-\alpha }f_{\alpha /2}(t^{-\alpha /2}r;n+1,-\frac{\alpha n}2-\alpha +1).
\end{multline*}

It is known (see Lemma 5 in \cite{Ps}) that
\begin{equation}
\left| f_{\alpha /2}(z;\mu ,\delta )\right| \le Cz^{1-\mu }\exp \{-\sigma z^{\frac2{2-\alpha}}\},
\end{equation}
if $\mu >1$. Therefore
$$
\left| \frac{d}{dr}\varphi_t(r)\right| \le Cr^{-n+1}t^{-\alpha}\exp \left\{ -\sigma \left( t^{-\alpha /2}r\right)^{\frac2{2-\alpha}}\right\}.
$$
Together with the inequalities (2.1)-(2.3), this implies (2.38). $\qquad \blacksquare$

\medskip
{\bf 2.4. Integral formulas.} It follows from Lemma 14 in \cite{Ps} that the following integral formulas hold:
\begin{equation}
\int\limits_{\Rn} Z_1^{(0)}(t,x-\xi ;\eta )\,dx=1;\quad \int\limits_{\Rn} Z_2^{(0)}(t,x-\xi ;\eta )\,dx=t;
\end{equation}
\begin{equation}
\int\limits_{\Rn} Y^{(0)}(t,x-\xi ;\eta )\,dx=\frac{t^{\alpha -1}}{\Gamma (\alpha )}.
\end{equation}

\medskip
{\it Remark}. The formula (2.46) is valid also for $0<\alpha <1$ providing a correction to the erroneous formula given in \cite{EK,EIK,K12}. This error does not influence other results from \cite{EK,K12}.

\section{Levi's Method}

{\bf 3.1. Construction.} We look for the functions $Z_1,Z_2,Y$ appearing in (1.6) assuming the integral representations:
\begin{equation}
Z_l(t,x;\xi )=Z_l^{(0)}(t,x-\xi ;\xi )+\int\limits_0^td\lambda
\int\limits_{\mathbb R^n}Y^{(0)}(t-\lambda ,x-y;y)Q_l(\lambda ,y;\xi
)\,dy;
\end{equation}
\begin{equation}
Q_l(t,x;\xi )=M_l(t,x;\xi )+\int\limits_0^td\lambda
\int\limits_{\mathbb R^n}K(t-\lambda,x;y)Q_l(\lambda ,y;\xi
)\,dy,
\end{equation}
where $l=1,2$,
\begin{multline}
M_l(t,x;\xi )=\sum\limits_{i,j=1}^n\left\{ [a_{ij}(x)-a_{ij}(\xi
)]\frac{\partial^2}{\partial x_i\partial x_j}Z_l^{(0)}(t,x-\xi ;\xi
)\right\}\\
+\sum\limits_{j=1}^nb_j(x)\frac{\partial Z_l^{(0)}(t,x-\xi ;\xi
)}{\partial x_j}+c(x)Z_l^{(0)}(t,x-\xi ;\xi );
\end{multline}
\begin{multline}
K(t,x;\xi )=\sum\limits_{i,j=1}^n\left\{ [a_{ij}(x)-a_{ij}(\xi
)]\frac{\partial^2}{\partial x_i\partial x_j}Y^{(0)}(t,x-\xi ;\xi
)\right\}\\
+\sum\limits_{j=1}^nb_j(x)\frac{\partial Y^{(0)}(t,x-\xi ;\xi
)}{\partial x_j}+c(x)Y^{(0)}(t,x-\xi ;\xi).
\end{multline}
Thus, the desired kernels consist of those for the equation with coefficients ``frozen'' at the parametric point $\xi$, plus correction terms constructed (as solutions of appropriate integral equations) in such a way that $LZ_l=0$ for $x\ne \xi$.

Similarly,
\begin{equation}
Y(t,x;\xi )=Y^{(0)}(t,x-\xi ;\xi )+\int\limits_0^td\lambda
\int\limits_{\mathbb R^n}Y^{(0)}(t-\lambda ,x-y;y)\Psi (\lambda ,y;\xi
)\,dy;
\end{equation}
\begin{equation}
\Psi (t,x;\xi )=K(t,x;\xi )+\int\limits_0^td\lambda
\int\limits_{\mathbb R^n}K(t-\lambda,x;y)\Psi (\lambda ,y;\xi )\,dy,
\end{equation}
so that $LY=0$ for $x\ne \xi$.

Let us consider first the case $n\ge 3$. It follows from (2.10)-(2.12) and (2.14) that
\begin{equation}
|M_1(t,x;\xi )|\le Ct^{-\alpha }|x-\xi |^{-n+\gamma }\rho_\sigma (t,x,\xi );
\end{equation}
\begin{equation}
|M_2(t,x;\xi )|\le Ct^{-\alpha +1}|x-\xi |^{-n+\gamma }\rho_\sigma (t,x,\xi );
\end{equation}
\begin{equation}
|K(t,x;\xi )|\le Ct^{-1}|x-\xi |^{-n+\gamma }\rho_\sigma (t,x,\xi ).
\end{equation}

We need to transform the estimates (3.7) and (3.9) in such a way that powers of $t$ become $>-1$, while powers of $|x-\xi |$ remain $>-n$. Of course, this is achieved at the expense of changing $\sigma$. We proceed as follows.

Since $\alpha <2$, we have $\alpha -1<\dfrac{\alpha}2$, so that there exists a constant $\nu_1\in (0,1)$, such that $\alpha -1<\dfrac{\nu_1\alpha}2$ or, equivalently, $\nu_1>2-\dfrac2\alpha$. Recalling our assumption (1.5), we can assume that
$$
\gamma > \nu_1>2-\frac2\alpha.
$$
In particular, $\dfrac{\nu_1\alpha}2-\alpha +1>0$, and we can define $\nu_0>0$ setting
$$
\nu_0=\frac{\frac12\nu_1\alpha -\alpha +1}{\frac12\alpha}=\nu_1-2+\frac2\alpha .
$$
Since $\nu_1<\gamma$ and $\alpha >1$, we have $\nu_0<\nu_1<\gamma$. Thus,
\begin{equation}
\frac{\nu_1\alpha}2-\alpha +1=\frac{\nu_0\alpha}2,\quad 0<\nu_0<\nu_1<\gamma .
\end{equation}

Now we can write
$$
t^{-\alpha}=t^{\frac{\nu_0\alpha}2-1}\left( t^{-\alpha /2}|x-\xi |\right)^{\nu_1}|x-\xi |^{-\nu_1},
$$
and decreasing $\sigma$ (we preserve the same letter) we obtain the estimate
\begin{equation}
|M_1(t,x;\xi )|\le Ct^{\frac{\nu_0\alpha}2-1}|x-\xi |^{-n+\gamma -\nu_1}\rho_\sigma (t,x,\xi );
\end{equation}
Similarly,
\begin{equation}
|K(t,x;\xi )|\le Ct^{\frac{\nu_0\alpha}2-1}|x-\xi |^{-n+\gamma -\nu_1}\rho_\sigma (t,x,\xi ).
\end{equation}

Now we can apply (without notable changes) the techniques from \cite{EK} (see also \cite{EIK}) to prove the solvability of the integral equations (3.2) and (3.6), and the estimates
\begin{equation}
|Q_1(t,x;\xi )|\le Ct^{\frac{\nu_0\alpha}2-1}|x-\xi |^{-n+\gamma -\nu_1}\rho_\sigma (t,x,\xi );
\end{equation}
\begin{equation}
|Q_2(t,x;\xi )|\le Ct^{-\alpha +1}|x-\xi |^{-n+\gamma -\nu_1}\rho_\sigma (t,x,\xi );
\end{equation}
\begin{equation}
|\Psi (t,x;\xi )|\le Ct^{\frac{\nu_0\alpha}2-1}|x-\xi |^{-n+\gamma -\nu_1}\rho_\sigma (t,x,\xi );
\end{equation}

Note that the techniques from \cite{EK,EIK} to prove (3.13)--(3.16) is more complicated than the standard method for second order parabolic equations. The estimates of iterated kernels are performed in two stages, like in the Levi method for parabolic systems \cite{E,F}.

For $n=2$, the above arguments carry over, if we roughen the estimates (2.21), (2.22), (2.24)--(2.26) substituting powers with the exponent $-\varepsilon$ with small $\varepsilon >0$ for the logarithmic factors. This leads to the same estimates (3.13)--(3.15) with smaller $\gamma$.

For $n=1$, we have in the same spirit that
\begin{equation}
|M_1(t,x;\xi )|\le Ct^{-(3-\gamma )\alpha /2}\rho_\sigma (t,x,\xi );
\end{equation}
\begin{equation}
|M_2(t,x;\xi )|\le Ct^{-(3-\gamma )\alpha /2+1}\rho_\sigma (t,x,\xi );
\end{equation}
\begin{equation}
|K(t,x;\xi )|\le Ct^{-(1-\gamma )\alpha /2-1}\rho_\sigma (t,x,\xi ).
\end{equation}
\begin{equation}
|\Psi (t,x;\xi )|\le Ct^{-(1-\gamma )\alpha /2-1}\rho_\sigma (t,x,\xi ).
\end{equation}
\begin{equation}
|Q_1(t,x;\xi )|\le Ct^{-(3-\gamma )\alpha /2}\rho_\sigma (t,x,\xi );
\end{equation}
\begin{equation}
|Q_2(t,x;\xi )|\le Ct^{-(3-\gamma )\alpha /2+1}\rho_\sigma (t,x,\xi );
\end{equation}
Note that the estimates (3.16), (3.18)--(3.20) have the same form as their counterparts for $0<\alpha <1$.

\medskip
{\bf 3.2. Increments.} The estimates of increments $\Delta_xM_j(t,x;\xi )=M_j(t,x;\xi )-M_j(t,x';\xi )$, $\Delta_xQ_j$ ($j=1,2$), $\Delta_xK$, $\Delta_x\Psi$ (the notation for increments should not be confused with that of the Laplacian) are derived just as in \cite{EK,EIK}. The results are as follows.

Let $x''$ be one of the points $x,x'$, for which $|x''-\xi |=\min (|x-\xi |,|x'-\xi |)$. If $n\ge 2$, then for some $\varepsilon >0$,
\begin{equation}
|\Delta_xM_1(t,x;\xi )|\le Ct^{-\alpha }|x-x'|^{\gamma -\varepsilon}|x''-\xi |^{-n+\varepsilon }\rho_\sigma (t,x'',\xi );
\end{equation}
\begin{equation}
|\Delta_xM_2(t,x;\xi )|\le Ct^{-\alpha +1}|x-x'|^{\gamma -\varepsilon}|x''-\xi |^{-n+\varepsilon }\rho_\sigma (t,x'',\xi );
\end{equation}
\begin{equation}
|\Delta_xK(t,x;\xi )|\le Ct^{-1}|x-x'|^{\gamma -\varepsilon}|x''-\xi |^{-n+\varepsilon }\rho_\sigma (t,x'',\xi );
\end{equation}
\begin{equation}
|\Delta_xQ_1(t,x;\xi )|\le Ct^{\frac{\nu_0\alpha}2-1}|x-x'|^{\gamma -\varepsilon}|x''-\xi |^{-n+\varepsilon }\rho_\sigma (t,x'',\xi );
\end{equation}
\begin{equation}
|\Delta_xQ_2(t,x;\xi )|\le Ct^{-\alpha +1}|x-x'|^{\gamma -\varepsilon}|x''-\xi |^{-n+\varepsilon }\rho_\sigma (t,x'',\xi );
\end{equation}
\begin{equation}
|\Delta_x\Psi (t,x;\xi )|\le Ct^{\frac{\nu_0\alpha}2-1}|x-x'|^{\gamma -\varepsilon}|x''-\xi |^{-n+\varepsilon }\rho_\sigma (t,x'',\xi ).
\end{equation}

If $n=1$, then for some $\varepsilon >0$,
\begin{equation}
|\Delta_xM_1(t,x;\xi )|\le Ct^{-\frac{3\alpha}2+\varepsilon }|x-x'|^{\varepsilon}\rho_\sigma (t,x'',\xi );
\end{equation}
\begin{equation}
|\Delta_xM_2(t,x;\xi )|\le Ct^{-\frac{3\alpha}2+\varepsilon+1 }|x-x'|^{\varepsilon}\rho_\sigma (t,x'',\xi );
\end{equation}
\begin{equation}
|\Delta_xK(t,x;\xi )|\le Ct^{-1-(1-\varepsilon )\frac{\alpha}2}|x-x'|^{\varepsilon}\rho_\sigma (t,x'',\xi );
\end{equation}
\begin{equation}
|\Delta_xQ_1(t,x;\xi )|\le Ct^{-\frac{3\alpha}2+\varepsilon }|x-x'|^{\varepsilon}\rho_\sigma (t,x'',\xi );
\end{equation}
\begin{equation}
|\Delta_xQ_2(t,x;\xi )|\le Ct^{-\frac{3\alpha}2+\varepsilon +1}|x-x'|^{\varepsilon}\rho_\sigma (t,x'',\xi );
\end{equation}
\begin{equation}
|\Delta_x\Psi (t,x;\xi )|\le Ct^{-\frac{\alpha}2-1+\varepsilon }|x-x'|^{\varepsilon}\rho_\sigma (t,x'',\xi ).
\end{equation}

In both cases the estimates for $\Delta_xM_1$, $\Delta_xM_2$, and $\Delta_xK$, that is (3.24)--(3.26) and (3.28)--(3.30) respectively, are obtained by elementary inequalities, and then applied to obtain the above estimates for $\Delta_xQ_1$, $\Delta_xQ_2$, and $\Delta_x\Psi$ using the general inequalities for convolution type integrals (see Sections 1.3.3 and 1.3.5 in \cite{EIK}).

The estimates of increments are important for the Levi method, since they make it possible to apply higher derivatives to expressions involving the fundamental solution.

\medskip
{\bf 3.3. Estimates for the fundamental solution.} The fact that the kernels we have constructed form indeed a fundamental solution will be proved later. Now we summarize their estimates.

\medskip
\begin{teo}
The Levi method kernels have the form
$$
Z_j(t,x,\xi )=Z_j^{(0)}(t,x-\xi ;\xi )+V_{Z_j}(t,x;\xi ),\quad j=1,2;
$$
$$
Y(t,x,\xi )=Y^{(0)}(t,x-\xi ;\xi )+V_Y(t,x;\xi ),
$$
where $Z_j^{(0)}$ and $Y^{(0)}$ satisfy the estimates (2.10)--(2.37). If $n\ge 2$, then
\begin{equation}
\left| D_x^m V_{Z_1}(t,x;\xi )\right| \le Ct^{\nu_0\alpha -1}|x-\xi |^{-n-|m|+(\gamma -\nu_1)+(2-\nu_0)}\rho_\sigma (t,x,\xi );
\end{equation}
\begin{equation}
\left| D_x^m V_{Z_2}(t,x;\xi )\right| \le Ct^{\frac{\nu_0\alpha}2+1-\alpha }|x-\xi |^{-n-|m|+(\gamma -\nu_1)+(2-\nu_0)}\rho_\sigma (t,x,\xi ),
\end{equation}
$|m|=0,1$;
\begin{equation}
\left| D_x^m V_{Z_1}(t,x;\xi )\right| \le Ct^{\gamma_1-\alpha }|x-\xi |^{-n+\gamma_2}\rho_\sigma (t,x,\xi );
\end{equation}
\begin{equation}
\left| D_x^m V_{Z_2}(t,x;\xi )\right| \le Ct^{\gamma_1-\alpha +1}|x-\xi |^{-n+\gamma_2}\rho_\sigma (t,x,\xi ),
\end{equation}
if $|m|=2$; here $\gamma_1$ and $\gamma_2$ are some positive constants;
\begin{equation}
\left| V_Y(t,x;\xi )\right| \le Ct^{\nu_0\alpha -1}|x-\xi |^{-n+(\gamma -\nu_1)+(2-\nu_0)}\rho_\sigma (t,x,\xi );
\end{equation}
\begin{equation}
\left| \frac{\partial }{\partial x_i}V_Y(t,x;\xi )\right| \le Ct^{\frac{\alpha \varkappa}2 +\frac{\nu_0\alpha}2 -1}|x-\xi |^{-n+1-\varkappa +\gamma -\nu_1}\rho_\sigma (t,x,\xi ),\quad n\ge 3,
\end{equation}
where $0<\varkappa <\gamma -\nu_1$. If $n=2$, then
\begin{equation}
\left| \frac{\partial }{\partial x_i}V_Y(t,x;\xi )\right| \le Ct^{\frac{\alpha \varkappa}2 +\frac{\nu_0\alpha}2 -1}|x-\xi |^{-1-\varkappa -\mu +\gamma -\nu_1}\rho_\sigma (t,x,\xi ),
\end{equation}
where $\mu >0$, $0<\varkappa +\mu <\gamma -\nu_1$;
\begin{equation}
\left| D_x^m V_Y(t,x;\xi )\right| \le Ct^{-1+\gamma_1}|x-\xi |^{-n+\gamma_2}\rho_\sigma (t,x,\xi ),
\end{equation}
if $|m|=2$; here $\gamma_1$ and $\gamma_2$ are positive constants.

If $n=1$, then
\begin{equation}
\left| D_x^m V_{Z_1}(t,x;\xi )\right| \le Ct^{(\gamma -m-1)\alpha/2}\rho_\sigma (t,x,\xi );
\end{equation}
\begin{equation}
\left| D_x^m V_{Z_2}(t,x;\xi )\right| \le Ct^{(\gamma -m-1)\alpha/2+1}\rho_\sigma (t,x,\xi );
\end{equation}
\begin{equation}
\left| D_x^m V_Y(t,x;\xi )\right| \le Ct^{\alpha -1+(\gamma -m-1)\alpha/2}\rho_\sigma (t,x,\xi ),
\end{equation}
m=0,1,2.
\end{teo}

\medskip
{\it Proof}. The ``additional'' terms $V_{Z_j},V_Y$ admit natural integral representations given by the second summands in the right-hand sides of (3.1) and (3.5). Estimates of the functions $Q_1,Q_2$ appearing in (3.1) and those for the function $\Psi$ from (3.5) are already available; see (3.13)--(3.15). Thus the estimates of $V_{Z_j}$ and their first derivatives follow from the estimates \cite{EIK} of special convolution operators; see Lemma 1.12 in \cite{EIK} for $n\ge 2$ and Lemma 1.5 in \cite{EIK} for $n=1$. These convolution operators contain singularities in the time and spatial variables, which may be rather strong; however they are compensated by the exponential factors $\rho_\sigma$. However, before using the above lemmas, the estimates must be prepared, ``taming'' the singularity. Let us show this procedure in the proof of the estimate (3.39) which will be used subsequently. To be specific, we assume that $n\ge 3$.

Our initial estimate of $\frac{\partial }{\partial x_i}Y^{(0)}$ is the inequality (2.13), while a bound for $\Psi$ is given by (3.15). Choose $\varkappa$ in such a way that $0<\varkappa <\gamma -\nu_1$ and write
\begin{multline*}
(t-\lambda )^{-\alpha -1}|x-y|^{-n+3}=\left[ (t-\lambda )^{-\alpha -\frac{\alpha \varkappa}2}|x-y|^{2+\varkappa}
\right] (t-\lambda )^{-1+\frac{\alpha \varkappa}2}|x-y|^{-n+1-\varkappa}\\
=\left[ (t-\lambda )^{-\alpha /2}|x-y|\right]^{2+\varkappa } (t-\lambda )^{-1+\frac{\alpha \varkappa}2}|x-y|^{-n+1-\varkappa}.
\end{multline*}
Changing $\sigma$ we can write
$$
\left| \frac{\partial}{\partial x_i}Y^{(0)}((t-\lambda ),x-y;y)\right| \le C(t-\lambda )^{-1+\frac{\alpha \varkappa}2}|x-y|^{-n+1-\varkappa}\rho_\sigma (t-\lambda,x,y).
$$
Now Lemma 1.12 from \cite{EIK} implies (3.39).

To prove the inequalities (3.36), (3.37), (3.41)--(3.44) for the second derivatives, one should repeat the reasoning from pages 344-345 in \cite{EIK} using the estimate (2.40) of the difference of values of $Y^{(0)}$ at different parameter points, as well as the estimates (3.25)--(3.27), (3.31)--(3.33) for increments of the functions $Q_1,Q_2,\Psi$. $\qquad \blacksquare$

\medskip
{\bf 3.4. Potentials.} For our situation, an analog of the heat potential is the function
\begin{equation}
W(t,x)=\int\limits_0^td\lambda \int\limits_{\Rn}Y^{(0)}(t-\lambda ,x-y;y)f(\lambda ,y)\,dy.
\end{equation}
We assume that $f(\lambda ,y)$ is a bounded function, jointly continuous in $(\lambda ,y)\in [0,T]\times \Rn$, and locally H\"older continuous in $y$, uniformly with respect to $\lambda$.

It is straightforward to check, using the estimates for $Y^{(0)}$, that the first derivatives in $x$ of $W(t,x)$ can be obtained by differentiating under the sign of integral in (3.45). Other derivatives are considered in the next proposition.

\medskip
\begin{prop}
The following differentiation formulas are valid:
\begin{equation}
\frac{\partial}{\partial t}W(t,x)=\int\limits_0^td\lambda \int\limits_{\Rn}\frac{\partial}{\partial t}Y^{(0)}(t-\lambda ,x-y;y)f(\lambda ,y)\,dy.
\end{equation}
\begin{multline}
\frac{\partial^2}{\partial x_i\partial x_j}W(t,x)=\int\limits_0^td\lambda \int\limits_{\Rn}\frac{\partial^2Y^{(0)}(t-\lambda ,x-y;y)}{\partial x_i\partial x_j}[f(\lambda ,y)-f(\lambda ,x)]\,dy\\
+\int\limits_0^tf(\lambda ,x)d\lambda \int\limits_{\Rn}\frac{\partial^2Y^{(0)}(t-\lambda ,x-y;y)}{\partial x_i\partial x_j}\,dy;
\end{multline}
\begin{multline}
\D W(t,x)=f(t,x)+\int\limits_0^td\lambda \int\limits_{\Rn}\frac{\partial Z_1^{(0)}(t-\lambda ,x-y;y)}{\partial t}[f(\lambda ,y)-f(\lambda ,x)]\,dy\\
+\int\limits_0^tf(\lambda ,x)d\lambda \int\limits_{\Rn}\frac{\partial Z_1^{(0)}(t-\lambda ,x-y;y)}{\partial t}\,dy.
\end{multline}
\end{prop}

\medskip
{\it Proof}. Denote
$$
W_h(t,x)=\int\limits_0^{t-h}d\lambda \int\limits_{\Rn}Y^{(0)}(t-\lambda ,x-y;y)f(\lambda ,y)\,dy.
$$
Then
\begin{multline*}
\frac{\partial}{\partial t}W_h(t,x)=\int\limits_{\Rn}Y^{(0)}(h,x-y;y)f(t-h,y)\,dy+\int\limits_0^{t-h}d\lambda \int\limits_{\Rn}\frac{\partial}{\partial t}Y^{(0)}(t-\lambda ,x-y;y)f(\lambda ,y)\,dy\\
\overset{\text{def}}{=}I_1+I_2.
\end{multline*}

Suppose, for example, that $n\ge 3$. By (2.14),
$$
|I_1|\le Ch^{-1}\int\limits_{\Rn}|y|^{-n+2}\exp \left\{ -\sigma |y|^{\frac2{2-\alpha}}t^{-\frac{\alpha}{2-\alpha}}\right\}\,dy \le Ch^{\alpha -1}\to 0,
$$
as $h\to 0$.

Let
$$
I_3=\int\limits_{\Rn}\frac{\partial}{\partial t}Y^{(0)}(t-\lambda ,x-y;y)f(\lambda ,y)\,dy.
$$
Using (2.17) we see that if $n=3$, then
$$
|I_3|\le C(t-\lambda )^{-\frac{\alpha}2-2}\int\limits_{\Rn}\exp \left\{ -\sigma (t^{-\alpha /2}|y|)^{\frac2{2-\alpha}}\right\} \,dy\le C(t-\lambda )^{\alpha -2}.
$$
The same estimate is obtained for $n=4$ and $n\ge 5$. Since $\alpha >1$, that means the convergence of $I_2$, as $h\to 0$, which implies (3.46).

It follows from (3.45) and (3.46) that
\begin{equation}
\lim\limits_{t\to 0}W(t,x)=\lim\limits_{t\to 0}\frac{\partial}{\partial t}W(t,x)=0.
\end{equation}
Then the relation between the Riemann-Liouville derivatives and the Laplace transform (see Theorem 7.3 in \cite{SKM}), together with the relation
$$
D_{0t}^{\alpha -1}Y^{(0)}(t,x-\xi ;\eta )=Z_1^{(0)}(t,x-\xi ;\eta )
$$
(see (2.6)), imply the identity $\D W=\dfrac{\partial v}{\partial t}$ where
$$
v(t,x)=\int\limits_0^td\lambda \int\limits_{\Rn}Z_1^{(0)}(t-\lambda ,x-y;y)f(\lambda ,y)\,dy.
$$

Now the rest of the proof of the equality (3.48) is identical to the one for $0<\alpha <1$ (see pages 247-248 in \cite{EK}, or pages 346-347 in \cite{EIK}).

The proof of (3.47) also repeats the reasoning in \cite{EK,EIK} (pages 243-246 and 342-345 respectively). $\qquad \blacksquare$

\medskip
Note that the formulas of Proposition 2 remain valid if the role of $f$ is played by the functions $Q_1,Q_2$, and $\Psi$. In these cases, instead of the simple H\"older condition, we use, in a similar way, the increment estimates (3.25)--(3.27), (3.31)--(3.33). The above results show that the Levi method constructions indeed produce solutions of the equation (1.3).

Note also that the differentiation formula (3.46) remains valid if
$$
|f(\lambda ,y)|\le C\lambda^{-a},\quad 0<a<1.
$$
This case will be important for checking the initial conditions in Section 4.

\section{Cauchy Problem}

{\bf 4.1. The case of zero initial functions.} Let us consider the equation (1.3), under the assumptions on coefficients formulated in the Introduction, with the initial conditions
\begin{equation}
u(0,x)=\frac{\partial u(0,x)}{\partial t}=0.
\end{equation}

\medskip
\begin{teo}
If $f$ is a bounded function, jointly continuous in $(t,x)$ and locally H\"older continuous in $x$, uniformly with respect to $t$, then
$$
u(t,x)=\int\limits_0^td\tau \int\limits_{\Rn}Y(t-\tau ,x;\xi )f(\tau ,\xi )\,d\xi
$$
is a bounded classical solution of the equation (1.3) satisfying the initial conditions (4.1).
\end{teo}

\medskip
{\it Proof}. The fact that $u$ is a solution follows from the Levi method construction. As for the initial conditions, we have the equalities (3.49); for the ``additional'' part
$$
u_{\text{add}}(t,x)=\int\limits_0^td\tau \int\limits_{\Rn}V_Y(t-\tau ,x;\xi )f(\tau ,\xi )\,d\xi ,
$$
we find, in a way similar to (3.46), that
$$
\frac{\partial u_{\text{add}}(t,x)}{\partial t}=\int\limits_0^td\tau \int\limits_{\Rn}\frac{\partial}{\partial t}V_Y(t-\tau ,x;\xi )f(\tau ,\xi )\,d\xi .
$$
Since $V_Y$ is less singular than $Y^{(0)}$, it is easy to check the initial conditions also for $u_{\text{add}}$, so that we obtain (4.1). $\qquad \blacksquare$

\medskip
{\bf 4.2. The homogeneous equation.} Here we consider the more complicated case of the equation (1.3) with $f=0$ and the initial conditions (1.2). We assume (in addition to the assumptions formulated in the Introduction) that:
\begin{description}
\item[(A)] $u_0(x)$ is bounded, continuously differentiable, and its first derivatives are bounded and
H\"older continuous with the exponent $\gamma_0>\dfrac{2-\alpha}\alpha$;
\item[(B)] $u_1(x)$ is bounded and continuous. If $n>1$, $u_1$ is H\"older continuous.
\item[(C)] The coefficients $a_{ij}$ are twice continuously differentiable with bounded derivatives of order
$\le 2$.
\end{description}

\bigskip
\begin{teo}
Under the above assumptions, the function
$$
u(t,x)=\int\limits_{\mathbb R^n} Z_1(t,x;\xi )u_0(\xi )\,d\xi +\int\limits_{\mathbb R^n} Z_2(t,x;\xi )u_1(\xi )\,d\xi
$$
is a bounded classical solution of the equation (1.3) (with $f=0$) satisfying the initial conditions (1.2).
\end{teo}

\medskip
{\it Proof}. Denote
$$
u^{(1)}(t,x)=\int\limits_{\mathbb R^n} Z_1(t,x;\xi )u_0(\xi )\,d\xi;
$$
$$
u^{(2)}(t,x)=\int\limits_{\mathbb R^n} Z_2(t,x;\xi )u_1(\xi )\,d\xi .
$$
Repeating the arguments from \cite{EIK} (pages 350-351) we show that $u^{(1)}(t,x)\to u_0(x)$, as $t\to 0$. Let us prove that $\dfrac{\partial u^{(1)}(t,x)}{\partial t}\to 0$, as $t\to 0$, for all $x\in \Rn$.

We have
\begin{multline*}
\dfrac{\partial u^{(1)}(t,x)}{\partial t}=\int\limits_{\mathbb R^n}\frac{\partial}{\partial t} Z_1^{(0)}(t,x-\xi ;\xi )u_0(\xi )\,d\xi +\frac{\partial}{\partial t}\int\limits_{\mathbb R^n}V_{Z_1}(t,x;\xi )u_0(\xi )\,d\xi \\
\overset{\text{def}}{=}u^{(1,1)}(t,x)+u^{(1,2)}(t,x).
\end{multline*}

Next,
\begin{multline*}
u^{(1,1)}(t,x)=\int\limits_{\mathbb R^n}\frac{\partial}{\partial t} Z_1^{(0)}(t,x-\xi ;x)u_0(\xi )\,d\xi \\
+\int\limits_{\mathbb R^n}\left[ \frac{\partial}{\partial t} Z_1^{(0)}(t,x-\xi ;\xi )
-\frac{\partial}{\partial t} Z_1^{(0)}(t,x-\xi ;x)\right] u_0(\xi )\,d\xi
\overset{\text{def}}{=}u^{(1,1,1)}(t,x)+u^{(1,1,2)}(t,x).
\end{multline*}
Pskhu \cite{Ps} proved, under our present assumptions regarding $u_0$, that $u^{(1,1,1)}(t,x)\to 0$, as $t\to 0$.

By (2.41),
\begin{equation}
\left| \frac{\partial}{\partial t} Z_1^{(0)}(t,x-\xi ;\xi )
-\frac{\partial}{\partial t} Z_1^{(0)}(t,x-\xi ;x)\right| \le Ct^{-\alpha -1}|x-\xi |^{-n+2+\gamma}\rho_\sigma (t,x,\xi ).
\end{equation}
By itself, this estimate is not sufficient for our purpose. It will work, if we obtain under the integral defining $u^{(1,1,2)}$ the difference $u_0(\xi )-u_0(x)$. To achieve that, we need an estimate for
$$
\frac{\partial}{\partial t}\int\limits_{\mathbb R^n}Z_1^{(0)}(t,x-\xi ;\xi )\,d\xi .
$$

Let us write the Taylor expansion of the function $\eta \mapsto \frac{\partial}{\partial t} Z_1^{(0)}(t,x-\xi ;\eta )$ on a neighborhood of the point $\eta =x$:
\begin{multline}
\frac{\partial}{\partial t} Z_1^{(0)}(t,x-\xi ;\eta)=\frac{\partial}{\partial t} Z_1^{(0)}(t,x-\xi ;x)+
(\eta -x)\cdot \nabla_\eta \frac{\partial}{\partial t} Z_1^{(0)}(t,x-\xi ;\eta )|_{\eta =x}\\
+\frac12\sum\limits_{j_1,j_2=1}^n\frac{\partial^2}{\partial \eta_{j_1}\partial \eta_{j_2}}\left[
\frac{\partial}{\partial t} Z_1^{(0)}(t,x-\xi ;\eta )\right]_{\eta =\eta_0}(\eta_{j_1}-x_{j_1})(\eta_{j_2}-x_{j_2})
\end{multline}
where $\eta_0$ lies on a segment joining $\eta$ and $x$.

In (4.3), we set $\eta =\xi$ and integrate in $\xi\in \Rn$. Then we notice that the integrals of the first two summands equal zero -- for the first one, it follows from (2.45) while the second is an odd function. Now we have to write the third summand more explicitly.

By (2.4), (2.6), (2.7), and (2.9) we find that
\begin{multline*}
\frac{\partial}{\partial t} Z_1^{(0)}(t,x-\xi ;\xi )=-\frac{\alpha n}2c_n\frac{t^{-\frac{\alpha n}2-1}}{(\det A(\eta ))^{1/2}}f_{\alpha /2}\left( (\mathfrak A(x-\xi ,\eta_0))^{1/2}t^{-\alpha /2};n-1,1-\frac{\alpha n}2\right) \\
+c_nt^{-\frac{\alpha n}2}(\det A(\eta ))^{-1/2}\frac{\partial}{\partial t}f_{\alpha /2}\left( (\mathfrak A(x-\xi ,\eta_0))^{1/2}t^{-\alpha /2};n-1,1-\frac{\alpha n}2\right) \\
\overset{\text{def}}{=}P_1(t,x-\xi ;\eta )+P_2(t,x-\xi ;\eta ),
\end{multline*}
and we have to substitute this into (4.3). We find using (2.8) that
\begin{multline*}
D^1_\eta f_{\alpha /2}\left( (\mathfrak A(x-\xi ,\eta_0))^{1/2}t^{-\alpha /2};n-1,1-\frac{\alpha n}2\right)
=f'_{\alpha /2}\left( (\mathfrak A(x-\xi ,\eta_0))^{1/2}t^{-\alpha /2};n-1,1-\frac{\alpha n}2\right) \\
\times \frac12 t^{-\alpha /2}(\mathfrak A(x-\xi ,\eta ))^{-1/2}\sum\limits_{i,j=1}^n D^1_\eta A^{(ij}(\eta ) (x_i-\xi_i)(x_j-\xi_j)\\
=-\frac14 t^{-\alpha }f_{\alpha /2}\left( (\mathfrak A(x-\xi ,\eta_0))^{1/2}t^{-\alpha /2};n+1,1-\frac{\alpha n}2-\alpha \right) \sum\limits_{i,j=1}^n D^1_\eta A^{(ij}(\eta ) (x_i-\xi_i)(x_j-\xi_j).
\end{multline*}

The next differentiation gives
\begin{multline*}
D^2_\eta f_{\alpha /2}\left( (\mathfrak A(x-\xi ,\eta_0))^{1/2}t^{-\alpha /2};n-1,1-\frac{\alpha n}2\right) \\
=-\frac14 t^{-\alpha }f_{\alpha /2}\left( (\mathfrak A(x-\xi ,\eta_0))^{1/2}t^{-\alpha /2};n+1,1-\frac{\alpha n}2-\alpha \right) \sum\limits_{i,j=1}^n D^2_\eta A^{(ij}(\eta ) (x_i-\xi_i)(x_j-\xi_j) \\
+\frac18 t^{-2\alpha }f_{\alpha /2}\left( (\mathfrak A(x-\xi ,\eta_0))^{1/2}t^{-\alpha /2};n+3,1-\frac{\alpha n}2-2\alpha \right) \left[ \sum\limits_{i,j=1}^n D^1_\eta A^{(ij}(\eta ) (x_i-\xi_i)(x_j-\xi_j)\right]^2
\end{multline*}
(in a certain abuse of notation we denote by $D^1_\eta$ and $D^2_\eta$ various partial derivatives; that does not influence the subsequent estimates).

Thus, we have found the summands of $D^2_\eta P_1$ (we should not forget other summands corresponding to differentiating $[\det A(\eta )]^{-1/2}$). It is known (see Lemma 5 in \cite{Ps}) that
\begin{equation}
|f_{\alpha /2}(z;\mu ,\delta )|\le C\exp \left\{ -\sigma z^{\frac2{2-\alpha}}\right\} \times \begin{cases}
1, &\text{if $0\le \mu <1$}\\
\log z, &\text{if $\mu =1$}\\
z^{1-\mu }, &\text{if $\mu >1$}.\end{cases}
\end{equation}
Using this inequality, we find that
\begin{equation}
\left| \left[ \sum\limits_{j_1,j_2=1}^n\frac{\partial^2}{\partial \eta_{j_1}\partial \eta_{j_2}}
P_1(t,x-\xi ;\eta_0)(\eta_{j_1}-x_{j_1})(\eta_{j_2}-x_{j_2})\right]_{\eta =\xi}\right|
\le Ct^{-\alpha -1}|x-\xi |^{-n+4}\rho_\sigma (t,x,\xi ),
\end{equation}
if $n\ge 2$ (the case $n=1$ leads to slightly different estimates with the same outcome).

Next, using (2.8) we find that
\begin{multline}
P_2(t,x-\xi ;\eta )\\
=\operatorname{const}\cdot t^{-\frac{\alpha n}2-\alpha -1}(\det A(\eta ))^{-1/2}\mathfrak A(x-\xi ,\eta)f_{\alpha /2}\left( (\mathfrak A(x-\xi ,\eta))^{1/2}t^{-\alpha /2};n+1,1-\frac{\alpha n}2-\alpha \right) .
\end{multline}
Differentiating the right-hand side of (4.6) twice in $\eta$ we obtain two kinds of terms:

1) The terms corresponding to differentiating $(\det A(\eta ))^{-1/2}\mathfrak A(x-\xi ,\eta)$. Their estimate is the same as the one for $P_2$, with the bound $Ct^{-\alpha -1}|x-\xi |^{-n+2}\rho_\sigma (t,x,\xi )$.

2) The terms corresponding to the first derivatives of the last factor in (4.6). Computing them as in the calculations with $P_1$ and using (4.4), we get the bound $Ct^{-2\alpha -1}|x-\xi |^{-n+4}\rho_\sigma (t,x,\xi )$ or, changing $\sigma$, the bound $Ct^{-\alpha -1}|x-\xi |^{-n+2}\rho_\sigma (t,x,\xi )$.

3) The terms corresponding to the second derivatives of the last factor in (4.6); they have the same bound as in the previous case.

Together with (4.5), these estimates show that
$$
\left| \int\limits_{\mathbb R^n}\frac{\partial}{\partial t}Z_1^{(0)}(t,x-\xi ;\xi )\,d\xi \right| \le Ct^{-\alpha -1}\int\limits_{\mathbb R^n}|x-\xi |^{-n+4}\rho_\sigma (t,x,\xi )\,d\xi ,
$$
so that
\begin{equation}
\left| \int\limits_{\mathbb R^n}\frac{\partial}{\partial t}Z_1^{(0)}(t,x-\xi ;\xi )\,d\xi \right| \le Ct^{\alpha -1}.
\end{equation}

Using (4.7) and the integral identity from (2.45) we obtain the asymptotic relation
\begin{equation}
u^{(1,1,2)}(t,x)=\int\limits_{\mathbb R^n}\left[ \frac{\partial}{\partial t} Z_1^{(0)}(t,x-\xi ;\xi )
-\frac{\partial}{\partial t} Z_1^{(0)}(t,x-\xi ;x)\right] [u_0(\xi )-u_0(x)]\,d\xi +o(1),\quad t\to 0.
\end{equation}

By our assumption (A),
$$
\left| u_0(\xi )-u_0(x)\right| \le C|x-\xi |.
$$
Using also the estimate (4.2), we find that
\begin{equation}
\left| u^{(1,1,2)}(t,x)\right| \le Ct^{-\alpha -1}\int\limits_{\mathbb R^n}|x-\xi |^{-n+3+\gamma }\rho_\sigma (t,x,\xi )\,d\xi +o(1).
\end{equation}
Under our assumption (C), we may set $\gamma =1$, and it follows from (4.9) that
$$
\left| u^{(1,1,2)}(t,x)\right| \le Ct^{\alpha -1}+o(1),
$$
so that
$u^{(1,1,2)}(t,x)\to 0$, as $t\to 0$, and we have proved that $u^{(1,1)}(t,x)\to 0$, as $t\to 0$.

Let us consider $u^{(1,2)}$. By definition of $V_{Z_1}$, we have
$$
\int\limits_{\mathbb R^n}V_{Z_1}(t,x;\xi )u_0(\xi )\,d\xi =\int\limits_0^td\lambda \int\limits_{\Rn}Y^{(0)}(t-\lambda ,x-y;y)F(\lambda ,y)\,dy
$$
where
\begin{equation}
F(\lambda ,y)=\int\limits_{\mathbb R^n}Q_1(\lambda ,y;\xi )u_0(\xi )\,d\xi .
\end{equation}

By (3.13), if $n\ge 2$, then
$$
|F(\lambda ,y)|\le C\lambda^{\frac{\nu_0\alpha }2-1}\int\limits_{\Rn}|y-\xi |^{-n+\gamma -\nu_1}\rho_\sigma (\lambda ,y,\xi )\,d\xi \le C\lambda^{\frac{\alpha \gamma}2+\frac{\alpha}2(\nu_0-\nu_1)-1}=C\lambda^{ \frac{\alpha \gamma}2-\alpha}
$$
where $\dfrac{\alpha \gamma}2-\alpha >-1$, and we may use the differentiation formula (3.46). We get
\begin{equation}
u^{(1,2)}(t,x)=\int\limits_0^td\lambda \int\limits_{\Rn}\frac{\partial}{\partial t}Y^{(0)}(t-\lambda ,x-y;y)F(\lambda ,y)\,dy.
\end{equation}

However, for studying the behavior of $u^{(1,2)}$, as $t\to 0$, the above estimate of $F$ is too rough, and that is not strange -- so far we have not used our assumptions (A) and (C). To deal with (4.10) in a more refined way, we need an estimate of the function
$$
q(\lambda ,y)=\int\limits_{\Rn}Q_1(\lambda ,y;\xi )\,d\xi .
$$
We assume that $\gamma =1$ and consider the case where $n\ge 3$; other cases are completely similar.

The first step is an estimate of the function
$$
I^{(1)}(\lambda ,y)=\int\limits_{\Rn}M_1(\lambda ,y;\xi )\,d\xi .
$$
We write $I^{(1)}$ as a sum of three integrals:
\begin{multline*}
I^{(1)}(\lambda ,y)=\int\limits_{\Rn}\left\{ \sum\limits_{i,j=1}^n[a_{ij}(x)-a_{ij}(\xi )]\left[
\frac{\partial^2}{\partial x_i\partial x_j}Z_1^{(0)}(t,x-\xi ;\xi )-\frac{\partial^2}{\partial x_i\partial x_j}Z_1^{(0)}(t,x-\xi ;x)\right]\right. \\
+\left. \sum\limits_{j=1}^nb_j(x)\left[
\frac{\partial}{\partial x_j}Z_1^{(0)}(t,x-\xi ;\xi )-\frac{\partial}{\partial x_j}Z_1^{(0)}(t,x-\xi ;x)\right]
+c(x)Z_1^{(0)}(t,x-\xi ;\xi )\right\} \,d\xi \\
+\int\limits_{\Rn}[a_{ij}(x)-a_{ij}(\xi )] \sum\limits_{i,j=1}^n\frac{\partial^2}{\partial x_i\partial x_j}Z_1^{(0)}(t,x-\xi ;x)\,d\xi \\
+\int\limits_{\Rn}b_j(x)\frac{\partial}{\partial x_j}Z_1^{(0)}(t,x-\xi ;x)\,d\xi \overset{\text{def}}{=}I^{(1,1)}+I^{(1,2)}+I^{(1,3)}.
\end{multline*}

Using (2.38) we find that
$$
\left| I^{(1,1)}\right| \le Ct^{-\alpha }\int\limits_{\Rn}|x-\xi |^{-n+2}\rho_\sigma (t,x,\xi )\,d\xi \le C.
$$

It follows from the explicit formulas for the equation (1.1) (see the proof of Lemma 15 in \cite{Ps}) that the function $\dfrac{\partial^2}{\partial y_i\partial y_j}Z_1^{(0)}(t,y;x)$ is even in $y$, while the function $\dfrac{\partial}{\partial y_j}Z_1^{(0)}(t,y;x)$ is odd in $y$. This means that $I^{(1,3)}=0$, whereas in $I^{(1,2)}$, after using the Taylor formula for $a_{ij}(\xi )$ in a neighborhood of $x$, the integrals with the first order terms are equal to zero. As a result,
$$
\left| I^{(1,2)}\right| \le Ct^{-\alpha }\int\limits_{\Rn}|x-\xi |^{-n+2}\rho_\sigma (t,x,\xi )\,d\xi \le C,
$$
so that
$\left| I^{(1)}(\lambda ,y)\right| \le C$.

For $q(\lambda ,y)$, integrating in $\xi$ the integral equation (3.2) of the Levi method we obtain the integral equation
$$
q(t,x)=I^{(1)}(t,x)+\int\limits_0^td\lambda \int\limits_{\Rn}K(t-\lambda ,x;y)q(\lambda ,y)\,dy.
$$
Solving it by iteration in a standard way, we obtain a representation
$$
q(t,x)=\sum\limits_{m=1}^\infty I^{(m)}(t,x)
$$
with
$$
\left| I^{(m)}(t,x)\right| \le C\frac{d^m}{\Gamma (\frac{m\alpha }2+1)}t^{(m-1)\alpha /2},\quad d>0.
$$
In particular, $|q(t,x)|\le C$ for $t\in [0,T],x\in \Rn$.

Returning to (4.10) we write
$$
F(\lambda ,y)=\int\limits_{\mathbb R^n}Q_1(\lambda ,y;\xi )[u_0(\xi )-u_0(x)]\,d\xi +u_0(x)q(\lambda ,y)
$$
and use the estimate (3.13) for $Q_1$ together with the inequality $|u_0(\xi )-u_0(x)|\le C|x-\xi |$. This shows that $|F(\lambda ,y)|\le C$.

Substituting this inequality into (4.11) and using (2.17) we find that
$$
\left| u^{(1,2)}(t,x)\right| \le C\int\limits_0^t(t-\lambda )^{\alpha -\frac{\alpha n}2-2}d\lambda \int\limits_{\Rn}\mu_n(t^{-\alpha /2}|x-y|)\rho_\sigma (t,x,y) \,dy\le Ct^{\alpha -1}\to 0,
$$
as $t\to 0$. Thus, we have proved that $\dfrac{\partial u_1(t,x)}{\partial t}\to 0$, as $t\to 0$, for all $x\in \Rn$.

Finally, we turn to $u^{(2)}(t,x)$ (again for $n\ge 3$). We have
$$
u^{(2)}(t,x)=u^{(2,1)}(t,x)+u^{(2,2)}(t,x)+u^{(2,3)}(t,x)
$$
where
$$
u^{(2,1)}(t,x)=\int\limits_{\mathbb R^n} Z_2^{(0)}(t,x-\xi ,x)u_1(\xi )\,d\xi ,
$$
$$
u^{(2,2)}(t,x)=\int\limits_{\mathbb R^n}\left[ Z_2^{(0)}(t,x-\xi ,\xi)-Z_2^{(0)}(t,x-\xi ,x)\right] u_1(\xi )\,d\xi ,
$$
$$
u^{(2,3)}(t,x)=\int\limits_0^t d\lambda \int\limits_{\mathbb R^n}Y^{(0)}(t-\lambda ,x-y;y)G(\lambda ,y)\,dy
$$
where
$$
G(\lambda ,y)=\int\limits_{\mathbb R^n}Q_2(\lambda ,y;\xi )u_1(\xi )\,d\xi .
$$

The investigation of $u^{(2,1)}$ is reduced by a change of variables to the case of the equation (1.1) studied by Pskhu \cite{Ps}. By his results,
$$
u^{(2,1)}(t,x)\longrightarrow 0,\quad \frac{\partial u^{(2,1)}(t,x)}{\partial t}\longrightarrow u_1(x),
$$
as $t\to 0$.

It follows from the definitions of $Z_1^{(0)}$ and $Z_2^{(0)}$ that
$$
\frac{\partial}{\partial t}\left[ Z_2^{(0)}(t,x-\xi ,\xi)-Z_2^{(0)}(t,x-\xi ,x)\right]=Z_1^{(0)}(t,x-\xi ,\xi)-Z_1^{(0)}(t,x-\xi ,x).
$$
Using (2.38) and (2.39) we see that
$$
\left| u^{(2,2)}(t,x)\right| \le Ct^{-\alpha +1}\int\limits_{\Rn}|x-\xi |^{-n+2+\gamma }\rho_\sigma (t,x,\xi)\,d\xi \le Ct^{1+\frac{\alpha \gamma}2}\to 0,\quad t\to 0;
$$
$$
\left| \frac{\partial u^{(2,2)}(t,x)}{\partial t}\right| \le Ct^{-\alpha }\int\limits_{\Rn}|x-\xi |^{-n+2+\gamma }\rho_\sigma (t,x,\xi)\,d\xi \le Ct^{\frac{\alpha \gamma}2}\to 0,\quad t\to 0;
$$

By (3.14), we have
$$
|G(\lambda ,y)| \le C\lambda^{1-\alpha }\int\limits_{\Rn}|y-\xi |^{-n+\gamma -\nu_1}\rho_\sigma (t,y,\xi)\,d\xi
\le C\lambda^{1-\alpha +\frac{\alpha}2(\gamma -\nu_1)}.
$$
Then it follows from (2.14) and (2.17) that
\begin{multline*}
\left| u^{(2,3)}(t,x)\right| \le C\int\limits_0^t \lambda^{1-\alpha +\frac{\alpha}2(\gamma -\nu_1)}(t-\lambda )^{-1}d\lambda \int\limits_{\mathbb R^n}|y-\xi |^{-n+2}\rho_\sigma (t,y,\xi)\,d\xi \\
\le C\int\limits_0^t \lambda^{1-\alpha +\frac{\alpha}2(\gamma -\nu_1)}(t-\lambda )^{-1+\alpha }d\lambda
\le Ct^{1+\frac{\alpha}2(\gamma -\nu_1)}\to 0;
\end{multline*}
\begin{multline*}
\left| \frac{\partial u^{(2,3)}(t,x)}{\partial t}\right| \le C\int\limits_0^t \lambda^{1-\alpha +\frac{\alpha}2(\gamma -\nu_1)}(t-\lambda )^{\alpha -\frac{\alpha n}2-2}d\lambda \int\limits_{\mathbb R^n}\mu_n( t^{-\alpha/2}|y-\xi |)\rho_\sigma (t,y,\xi)\,d\xi \\
\le C\int\limits_0^t \lambda^{1-\alpha +\frac{\alpha}2(\gamma -\nu_1)}(t-\lambda )^{\alpha -2}d\lambda \le Ct^{\frac{\alpha}2(\gamma -\nu_1)}\to 0,
\end{multline*}
as $t\to 0$. $\qquad \blacksquare$

\section{Uniqueness Theorem}

{\bf 5.1. The adjoint problem.} For equations with variable coefficients and $0<\alpha <1$, uniqueness theorems were proved in \cite{EK,EIK,K90} using the maximum principle arguments. For $1<\alpha <2$, the structure of the fractional derivative $\D$ is different -- its Marchaud forms (see \cite{SKM}) contain either the first derivative or the second difference, thus being not suitable for the maximum principle. Note also that the positivity of the function $Y^{(0)}$ is violated for $1<\alpha <2$, $n\ge 4$ \cite{Ps}.

Therefore it is natural to try for $1<\alpha <2$ another classical method \cite{F} based on the representation of solutions using a fundamental solution of the adjoint problem. Classically, the adjoint problem is a kind of the Cauchy problem with data on the right end of the interval. The time derivative is preserved in the adjoint operator, only with a different sign.

In the fractional case, it is known \cite{PC} that the operator, adjoint to the (left-sided) Caputo-Dzhrbashyan fractional derivative $\mathbb D^{(\alpha )}$ is the right-sided Riemann-Liouville fractional derivative. In contrast to the classical situation, the form of the latter operator depends on the interval on which the adjoint problem is considered. This makes the use of adjoints more complicated necessitating their subtler definition. Such a definition was proposed by Pskhu \cite{Ps} for the equation (1.1). Below we adapt his approach to our general case.

In this section we assume, in addition to the assumptions from Introduction, that the following holds:
\begin{description}
\item(D) All the functions
$$
\frac{\partial a_{ij}}{\partial x_k},\quad \frac{\partial^2 a_{ij}}{\partial x_k\partial x_l},\quad \frac{\partial b_i}{\partial x_k}\quad (i,j,k,l=1,\ldots ,n)
$$
exist, are bounded and uniformly H\"older continuous.
\end{description}

In the adjoint operator $\mathcal B^*$ with respect to the spatial variables,
$$
\mathcal B^* u(t,x)=\sum\limits_{i,j=1}^n a_{ij}(x)\frac{\partial^2 u(t,x)}{\partial x_i\partial x_j}+\sum\limits_{k=1}^nb_k^*(x)\frac{\partial u(t,x)}{\partial x_k}+c^*(x)u(t,x),
$$
the higher coefficients $a_{ij}$ are the same as in $\mathcal B$ (we have assumed that $a_{ij}=a_{ji}$),
$$
b_i^*(x)=-b_i(x)+2\sum\limits_{j=1}^n \frac{\partial a_{ij}(x)}{\partial x_j},
$$
$$
c^*(x)=c(x)-\sum\limits_{i=1}^n\frac{\partial b_i(x)}{\partial x_i}+\sum\limits_{i,j=1}^n \frac{\partial^2 a_{ij}(x)}{\partial x_i\partial x_j}.
$$

Let $S\subset \Rn$ be a bounded domain with a smooth boundary $\partial S$. Denote $E=\{ (t,x):\ 0<t<T,x\in S\}$, $E_t=\{ (\eta ,x):\ 0<\eta <t,x\in S\}$. In agreement with the earlier definition, we call $u(t,x)$ a classical solution of the equation (1.3) on $E$, if: 1) it satisfies on $\overline{S}$ the conditions (i)-(ii) from the definition of the classical solution from Introduction; in particular, for each $x\in \overline{S}$, there exist continuous on $\overline{S}$ limits $u_0,u_1$ of $u(t,x)$ and $\dfrac{\partial u(t,x)}{\partial t}$ respectively, as $t\to 0$; 2) $u(t,x)$ satisfies (1.3) at all the points $(t,x)\in E$.

The adjoint operator $L^*$ acts on functions $v(t,x;\eta ,\xi )$, in the variables $\eta ,\xi$, where $\eta <t$, as follows:
$$
L^*v=D_{t\eta }^\alpha v(t,x;\eta ,\xi )-\mathcal B^*_\xi v(t,x;\eta ,\xi ).
$$
Here $D_{t\eta }^\alpha$ (see the definition in Section 2.1 taken from \cite{Ps}) is the right-sided Riemann-Liouville derivative in the variable $\eta$ with the base point $t$.

If we consider the terminal value problem for the equation $L^*v=g$, $g=g(t,x;\eta ,\xi )$, with zero terminal condition at $\eta =t$, this problem is equivalent, via time reflection, to the homogeneous Cauchy problem considered above, with $\mathcal B^*$ substituted for $\mathcal B$. Therefore under our condition (D), for such a terminal value problem there exists a fundamental solution
\begin{equation}
Y_*(t,x;\eta ,\xi )=Y_*^{(0)}(t-\eta ,x-\xi ,x)+V_{Y_*}(t,x;\eta ,\xi )
\end{equation}
satisfying the same estimates as the function $Y$ above. Below, dealing with estimates for $Y_*^{(0)}$ and $V_{Y_*}$ we will refer to the appropriate estimates for $Y^{(0)}$ and $V_Y$. One should only remember that the operators will act on the function (5.1) in the variable $\xi$, and $x$ will be the integration variable.

Suppose that a function $v(t,x;\eta ,\xi )$ is continuous on $Q\times \overline{E_t}$, $Q\subset E$, together with its first and second derivatives in $\xi$ and its fractional derivatives $D_{t\eta }^\alpha$ and $D_{t\eta }^{\alpha -1}$. We also assume that
\begin{equation}
\lim\limits_{\eta \to t}\left( D_{t\eta }^{\alpha -1}v\right) (t,x;\eta ,\xi )=\lim\limits_{\eta \to t}\left( D_{t\eta }^{\alpha -2}v\right) (t,x;\eta ,\xi )=0.
\end{equation}
Denote
$$
w(t,x;\eta ,\xi )=v(t,x;\eta ,\xi )+Y_*(t,x;\eta ,\xi ).
$$

We will use the Green formula for the elliptic operator $\mathcal B$ (see \cite{Mi}). Let $\Omega$ be a smooth domain in $\Rn$. Denote by $X_i$ the direction cosines of the outer normal to $\Omega$. For $\xi \in \partial \Omega$, $\nu_\xi$ will denote the conormal at $\xi$, that is a vector with the direction cosines
$$
Y_i=\frac1{a(\xi )}\sum\limits_{k=1}^na_{ik}(\xi )X_k(\xi ),\quad a(\xi )=\left[ \sum\limits_{i=1}^n\left( \sum\limits_{k=1}^na_{ik}(\xi )X_k(\xi )\right)^2\right]^{1/2}.
$$

For smooth functions $U,V$ on $\overline{\Omega }$,
\begin{equation}
\int\limits_{\Omega }(V\cdot \mathcal BU-U\cdot \mathcal B^*V)\,d\Omega =\int\limits_{\partial \Omega }\left[ a\left( V\frac{\partial U}{\partial \nu}-U\frac{\partial V}{\partial \nu}\right) +bUV\right]\,dS_{\partial \Omega }
\end{equation}
where
$$
b(\xi )=\sum\limits_{i=1}^n\left( b_i(\xi )-\sum\limits_{k=1}^n\frac{\partial a_{ik}(\xi )}{\partial x_k}\right) X_i(\xi ),\quad \xi \in \partial \Omega.
$$

\medskip
\begin{prop}
Under the above assumptions, for each $(t,x)\in Q$, a classical solution $u$ of the equation (1.3) on $E$ has the representation
\begin{multline}
u(t,x)=\int\limits_S u_0(\xi )\left[ D_{t\eta }^{\alpha -1}w(t,x;\eta ,\xi )\right]_{\eta =0}\,d\xi +\int\limits_S u_1(\xi )\left[ D_{t\eta }^{\alpha -2}w(t,x;\eta ,\xi )\right]_{\eta =0}\,d\xi \\
+G(u;t,x)+F(u,f,g;t,x)
\end{multline}
where
\begin{multline*}
G(u;t,x)=\int\limits_0^td\eta \int\limits_{\partial S}\biggl\{ a(\xi )\left[ w(t,x;\eta ,\xi )\frac{\partial}{\partial \nu_\xi}u(t,\xi )-u(t,\xi )\frac{\partial}{\partial \nu_\xi}w(t,x;\eta ,\xi )\right]  \\
+b(\xi )w(t,x;\eta ,\xi )u(t,\xi )\biggr\}\,dS_\xi ,
\end{multline*}
$$
F(u,f,g;t,x)=\int\limits_0^td\eta \int\limits_S[w(t,x;\eta ,\xi )f(\eta ,\xi )-u(\eta ,\xi )g(t,x;\eta ,\xi )]\,d\xi.
$$
\end{prop}

\medskip
{\it Proof}. Let $S^\varepsilon \subset S$ ($\varepsilon >0$) be a smooth domain, such that $\operatorname{dist} (\partial S^{\varepsilon},\partial S)\le \varepsilon$. Denote $S_r^\varepsilon =S^\varepsilon \setminus B_r$, $B_r=\{ \xi \in S:\ |x-\xi |<r\}$. Let us consider the expression
\begin{equation}
\int\limits_\delta^td\eta \int\limits_{S_r^\varepsilon}[w(t,x;\eta ,\xi )(Lu)(\eta ,\xi )-u(\eta ,\xi )L^*w(t,x;\eta ,\xi )]\,d\xi
\end{equation}
with $0<\delta <t$, $\overline{B_r}\subset S^\varepsilon$.

Following \cite{Ps}, denote
$$
\left( I_{s\delta}^\zeta g\right) (\eta )=\begin{cases}
\dfrac{\operatorname{sign}(\delta -s)}{\Gamma (-\zeta )}\int\limits_s^\delta g(t)|\eta-t|^{-\zeta -1}dt, &\text{ if $\zeta <0$};\\
0, &\text{ if $\zeta =0$};\end{cases}
$$
$$
R_1=-\left[ \left( D_{t\eta }^{\alpha -1}w(t,x;\eta ,\xi )\right) u(\eta ,\xi )\right]_{\eta =\delta};
$$
$$
R_2=\int\limits_\delta^t w(t,x;\eta ,\xi )I_{0\delta}^{\alpha -2}\frac{\partial^2}{\partial \eta^2}u(\eta ,\xi )\,d\eta -\left[ \left( D_{t\eta }^{\alpha -2}w(t,x;\eta ,\xi )\right) \frac{\partial}{\partial \eta}u(\eta ,\xi )\right]_{\eta =\delta}.
$$
For each $\xi\in S_r^\varepsilon$,
\begin{multline*}
\int\limits_\delta^t w(t,x;\eta ,\xi )\left( \mathbb D_\eta^{(\alpha )}u\right) (\eta ,\xi )\,d\eta =\int\limits_\delta^t w(t,x;\eta ,\xi )\left[ I_{0\delta}^{\alpha -2}+D_{\delta \eta}^{\alpha -2}\right]\frac{\partial^2}{\partial \eta^2}u(\eta ,\xi )\,d\eta \\
=\int\limits_\delta^t \left( D_{t\eta}^{\alpha -2}w\right) (t,x;\eta ,\xi )\frac{\partial^2}{\partial \eta^2}u(\eta ,\xi )\,d\eta +\int\limits_\delta^t w(t,x;\eta ,\xi )I_{0\delta}^{\alpha -2}\frac{\partial^2}{\partial \eta^2}u(\eta ,\xi )\,d\eta \\
=\int\limits_\delta^t \left( D_{t\eta}^{\alpha -1}w\right) (t,x;\eta ,\xi )\frac{\partial}{\partial \eta}u(\eta ,\xi )\,d\eta +R_2.
\end{multline*}
We have used the classical and fractional versions of integration by parts \cite{SKM,KST}, the assumption (5.2), and the fact that $D_{t\eta}^{\alpha -2}Y_*(t,x;\eta ,\xi )\to 0$ for $\eta \to t$, if $\xi \in S_r^\varepsilon$ (so that $x-\xi $ is separated from zero).

Integrating by parts again we find that
$$
\int\limits_\delta^t w(t,x;\eta ,\xi )\left( \mathbb D_\eta^{(\alpha )}u\right) (\eta ,\xi )\,d\eta =\int\limits_\delta^t \left( D_{t\eta}^\alpha w(t,x;\eta ,\xi )\right) u(\eta ,\xi )\,d\eta +R_1+R_2.
$$
Now, using also the Green formula (5.3) we obtain that the expression (5.5) equals
\begin{multline}
\int\limits_{S_r^\varepsilon}(R_1+R_2)\,d\xi -\int\limits_\delta^t\int\limits_{\partial S^\varepsilon}\left[ a\left( w\frac{\partial u}{\partial \nu}-u\frac{\partial w}{\partial \nu}\right) +buw\right] \,dS \\
+\int\limits_\delta^t\int\limits_{\partial B_r}\left[ a\left( w\frac{\partial u}{\partial \nu}-u\frac{\partial w}{\partial \nu}\right) +buw\right] \,dS.
\end{multline}

Considering the integrals of $R_1$ and $R_2$ we note that the main singular part of $Y_*$ coincides, up to a change of variables, with the fundamental solution of the Cauchy problem for the model equation (1.1); this main part (actually, its fractional derivatives) is ``responsible'' for the solution of the Cauchy problem (see Theorem 3 and its proof). The additional term $V_{Y_*}$ is less singular. As a result, we can repeat the reasoning from \cite{Ps} to obtain the equality
\begin{equation}
\lim\limits_{\varepsilon \to 0}\lim\limits_{\delta \to 0}\lim\limits_{r \to 0}\int\limits_{S_r^\varepsilon}(R_1+R_2)\,d\xi =-\int\limits_Su_0(\xi )\left( D_{t\eta}^{\alpha -1}w(t,x;\eta ,\xi )\right) \,d\xi -\int\limits_Su_1(\xi )\left( D_{t\eta}^{\alpha -2}w(t,x;\eta ,\xi )\right) \,d\xi .
\end{equation}

Due to the smoothness of the functions $u$ and $v$,
\begin{multline*}
\lim\limits_{r \to 0}\int\limits_\delta^t d\eta \int\limits_{\partial B_r}a\left( w\frac{\partial u}{\partial \nu}-u\frac{\partial w}{\partial \nu}\right) \,dS \\
=\lim\limits_{r \to 0}\int\limits_\delta^t d\eta \int\limits_{\partial B_r}a(\xi )\left[ Y_*(t,x;\eta ,\xi )\frac{\partial u(\eta ,\xi )}{\partial \nu_\xi }-u(\eta ,\xi )\frac{\partial }{\partial \nu_\xi }Y_*(t,x;\eta ,\xi )\right]\,dS_\xi .
\end{multline*}

It follows directly from the estimates (2.12), (2.23) and (2.33) for $Y^{(0)}$ and the estimates (3.38), (3.44) for $V_Y$ that
$$
\lim\limits_{r \to 0}\int\limits_\delta^t d\eta \int\limits_{\partial B_r}a(\xi )Y_*(t,x;\eta ,\xi )\frac{\partial u(\eta ,\xi )}{\partial \nu_\xi }\,dS_\xi =0.
$$
Similarly, by (2.13), (2.24), (2.33), (3.39), and (3.44), for an arbitrary $t^*<t$,
$$
\lim\limits_{r \to 0}\int\limits_\delta^{t^*} d\eta \int\limits_{\partial B_r}a(\xi )u(\eta ,\xi )\frac{\partial }{\partial \nu_\xi }Y_*(t,x;\eta ,\xi )\,dS_\xi =0.
$$
Therefore
\begin{multline*}
\lim\limits_{r \to 0}\int\limits_\delta^t d\eta \int\limits_{\partial B_r}a(\xi )u(\eta ,\xi )\frac{\partial }{\partial \nu_\xi }Y_*(t,x;\eta ,\xi )\,dS_\xi \\
=\lim\limits_{r \to 0}\int\limits_{t^*}^td\eta \int\limits_{\partial B_r}a(\xi )[u(\eta ,\xi )-u(t,x)]\frac{\partial }{\partial \nu_\xi }Y_*(t,x;\eta ,\xi )\,dS_\xi \\
+u(t,x)\lim\limits_{r \to 0}\int\limits_{t^*}^t d\eta \int\limits_{\partial B_r}a(\xi )\frac{\partial }{\partial \nu_\xi }Y_*(t,x;\eta ,\xi )\,dS_\xi \overset{\text{def}}{=}\lim\limits_{r \to 0}I_1+u(t,x)\lim\limits_{r \to 0}I_2.
\end{multline*}

Let us study $I_1$ taking into account that
$$
|u(\eta ,\xi )-u(t,x)|\le C(|t-\eta |+|x-\xi |).
$$
It is important here to use the precise estimate (2.13); the estimate (2.14) used in the Levi method is not sufficient. We have for $n\ge 3$, by (2.13) and (3.39),
\begin{multline*}
|I_1|\le C\int\limits_{t^*}^t (t-\eta )^{-\alpha }d\eta \int\limits_{\partial B_r}|x-\xi |^{-n+3}\rho_\sigma (t-\eta ,x,\xi )\,dS_\xi \\
+C\int\limits_{t^*}^t (t-\eta )^{\frac{\alpha \varkappa}2+\frac{\nu_0\alpha}2}d\eta \int\limits_{\partial B_r}|x-\xi |^{-n+1+\gamma -\nu_1-\varkappa}\rho_\sigma (t-\eta ,x,\xi )\,dS_\xi \\
+C\int\limits_{t^*}^t (t-\eta )^{-1+\frac{\alpha \varkappa}2}d\eta \int\limits_{\partial B_r}|x-\xi |^{-n+2-\varkappa}\rho_\sigma (t-\eta ,x,\xi )\,dS_\xi \\
+C\int\limits_{t^*}^t (t-\eta )^{\frac{\alpha \varkappa}2+\frac{\nu_0\alpha}2-1}d\eta \int\limits_{\partial B_r}|x-\xi |^{-n+2+\gamma -\nu_1-\varkappa}\rho_\sigma (t-\eta ,x,\xi )\,dS_\xi .
\end{multline*}

In each of the integrals, $|x-\xi |=r$, so that they are easily calculated showing that $I_1\to 0$, as $r\to 0$. The case $n=2$ is similar. If $n=1$, then
\begin{multline*}
I_1=\int\limits_{t^*}^t\left\{ a(x+r)[u(\eta ,x+r)-u(t,x)]\frac{\partial }{\partial r}Y_*(t,x;\eta ,x+r)\right. \\
\left. -a(x-r)[u(\eta ,x-r)-u(t,x)]\frac{\partial }{\partial r}Y_*(t,x;\eta ,x-r)\right\} \,d\eta .
\end{multline*}
We can write $I_1=I_1^{(1)}+I_1^{(2)}$ where the summands correspond to the decomposition $Y_*=Y_*^{(0)}+V_{Y_*}$. Note that $Y_*^{(0)}(t-\eta ,r;x)=Y_*^{(0)}(t-\eta ,-r;x)$. Therefore
\begin{multline*}
I_1^{(1)}=\int\limits_{t^*}^t [a(x+r)u(\eta ,x+r)-a(x-r)u(\eta ,x-r)]\frac{\partial }{\partial r}Y_*^{(0)}(t-\eta ,r;x)\,d\eta \\
+u(t,x)\int\limits_{t^*}^t [a(x-r)-a(x+r)]\frac{\partial }{\partial r}Y_*^{(0)}(t-\eta ,r;x)\,d\eta .
\end{multline*}

Using (2.34) we find that
$$
\left| I_1^{(1)}\right| \le Cr^2\int\limits_{t^*}^t(t-\eta )^{-\frac{\alpha}2-1}\rho_\sigma (t-\eta ,r,0)\,d\eta .
$$
Taking $\varkappa \in (0,1)$ we write
$$
(t-\eta )^{-\frac{\alpha}2-1}r^2=\left[ (t-\eta )^{-\frac{\alpha}2}r\right]^{2-\varkappa }(t-\eta )^{\frac{\alpha}2-1-\frac{\varkappa \alpha}2}r^\varkappa .
$$
Changing $\sigma$ and removing the first factor, and then bounding $\rho_\sigma$ by 1, we obtain that
$$
\left| I_1^{(1)}\right| \le Cr^\varkappa \int\limits_{t^*}^t(t-\eta )^{(1-\varkappa )\frac{\alpha}2-1}d\eta \to 0,
$$
as $r\to 0$.

The fact that $I_1^{(2)}\to 0$, as $r\to 0$, follows from the continuity of $u$, the estimate (3.44) for the first derivative of $V_{Y_*}$, and the dominated convergence theorem.

Turning to $I_2$ we write
\begin{multline*}
I_2=\int\limits_{t^*}^t d\eta \int\limits_{\partial B_r}a(\xi )\frac{\partial }{\partial \nu_\xi }Y_*(t,x;\eta ,x)\,dS_\xi +\int\limits_{t^*}^t d\eta\int\limits_{\partial B_r}a(\xi )\left[ \frac{\partial }{\partial \nu }Y_*(t,x;\eta ,\xi )-\right.\\
\left. -\frac{\partial }{\partial \nu }Y_*(t,x;\eta ,x)\right] \,dS_\xi +\int\limits_{t^*}^t d\eta \int\limits_{\partial B_r}a(\xi )\frac{\partial }{\partial \nu_\xi }V_{Y_*}(t,x;\eta ,\xi )\,dS_\xi \overset{\text{def}}{=}I_2^{(1)}+I_2^{(2)}+I_2^{(3)}.
\end{multline*}

Using the estimates (2.40), (3.39), and (3.43) we show in a straightforward way that $I_2^{(2)}\to 0$ and $I_2^{(3)}\to 0$, as $r\to 0$.

Denote by $\nu_\xi^{(x)}$ the conormal at a point $\xi$ corresponding to the operator $\mathcal B^*_\xi (x)$ acting in $\xi$, with the leading coefficients ``frozen'' at the point $x$ and other terms discarded. Then
\begin{multline*}
I_2^{(1)}=\int\limits_{t^*}^t d\eta \int\limits_{\partial B_r}a(\xi )\frac{\partial }{\partial \nu_\xi^{(x)} }Y_*^{(0)}(t-\eta ,x-\xi ;x)\,dS_\xi +\int\limits_{t^*}^t d\eta \int\limits_{\partial B_r}a(\xi )\left[ \frac{\partial }{\partial \nu_\xi} Y_*^{(0)}(t-\eta ,x-\xi ;x)\right. \\
\left. -\frac{\partial }{\partial \nu_\xi^{(x)}}Y_*^{(0)}(t-\eta ,x-\xi ;x)\right] \,dS_\xi \overset{\text{def}}{=}I_2^{(1,1)}+I_2^{(1,2)}.
\end{multline*}

The second term is simpler -- for $n\ge3$ (other cases are similar), recalling the definition of a conormal and transforming the estimate (2.13) as we did treating $I_1$ above, we find that, with $0<\varkappa <\gamma$,
$$
\left| I_2^{(1,2)}\right| \le C\int\limits_{t^*}^t (t-\eta )^{-1+\frac{\alpha \varkappa}2}d\eta \int\limits_{\partial B_r}|x-\xi |^{-n+1-\varkappa +\gamma }\rho_\sigma (t-\eta ,x,\xi )\,dS_\xi \le Cr^{\gamma -\varkappa}\to 0,
$$
as $r\to 0$.

Next,
\begin{multline*}
I_2^{(1,1)}=\int\limits_{t^*}^t d\eta \int\limits_{\partial B_r}[a(\xi )-a(x)]\frac{\partial }{\partial \nu_\xi^{(x)} }Y_*^{(0)}(t-\eta ,x-\xi ;x)\,dS_\xi \\
+a(x)\int\limits_{t^*}^t d\eta \int\limits_{\partial B_r}\frac{\partial }{\partial \nu_\xi^{(x)} }Y_*^{(0)}(t-\eta ,x-\xi ;x)\,dS_\xi
\overset{\text{def}}{=}I_2^{(1,1,1)}+I_2^{(1,1,2)}.
\end{multline*}
As before, an easy estimate shows that $I_2^{(1,1,1)}\to 0$, as $r\to 0$. In the last term we use the Green formula (actually, for the operator with constant coefficients). This results in the representation
$$
I_2^{(1,1,2)}=-\int\limits_{t^*}^t d\eta \int\limits_{B_R\setminus B_r}\mathcal B^*_\xi (x)Y_*^{(0)}(t-\eta ,x-\xi ;x)\,d\xi +a(x)\int\limits_{t^*}^t d\eta \int\limits_{\partial B_R}\frac{\partial }{\partial \nu_\xi^{(x)} }Y_*^{(0)}(t-\eta ,x-\xi ;x)\,dS_\xi .
$$
Here $R>r$, and we will pass to the limit, as $R\to \infty$ and $r\to 0$. We have
\begin{multline*}
\left| a(x)\int\limits_{t^*}^t d\eta \int\limits_{\partial B_R}\frac{\partial }{\partial \nu_\xi^{(x)} }Y_*^{(0)}(t-\eta ,x-\xi ;x)\,dS_\xi \right| \\
\le C\int\limits_{t^*}^t (t-\eta )^{-1}d\eta \int\limits_{\partial B_R}|x-\xi |^{-n+1}\rho_\sigma (t-\eta ,x,\xi )\,dS_\xi \\
=\operatorname{const}\cdot \int\limits_{t^*}^t (t-\eta )^{-1}\exp \left\{ -\sigma \left[ (t-\eta )^{-\alpha /2}R\right]^{\frac2{2-\alpha}}\right\}\,d\eta \\
=\operatorname{const}\cdot \int\limits_0^{(t-t^*)R^{-2/\alpha}}\tau^{-1}\exp \left\{ -\sigma \tau^{-\frac{\alpha}{2-\alpha}}\right\} \,d\tau \to 0,
\end{multline*}
as $R\to \infty$. Finally, by (2.6) and (2.45),
\begin{multline*}
\int\limits_{t^*}^t d\eta \int\limits_{B_R\setminus B_r}\mathcal B^*_\xi (x)Y_*^{(0)}(t-\eta ,x-\xi ;x)\,d\xi =\int\limits_{B_R\setminus B_r}D_{t\eta}^{\alpha -1}Y^{(0)}(t-\eta ,x-\xi ;x)|_{\eta =t^*}\,d\xi \\
=\int\limits_{B_R\setminus B_r}Z_1^{(0)}(t-t^*,x-\xi ;x)\,d\xi \to 1,
\end{multline*}
as $R\to \infty$, $r\to 0$.

Hence, $I_2^{(1,1)}\to -1$, as $r\to 0$, so that $I_2\to -1$, as $r\to 0$. Therefore the identity of (5.5) and (5.6), together with (5.7) and the above limit relations, results in (5.4). $\qquad \blacksquare$

\bigskip
{\bf 5.2. Uniqueness conditions.} Using Proposition 3 we can repeat the arguments from \cite{Ps} and obtain the following uniqueness result.

\medskip
\begin{teo}
Suppose that the condition (D) is satisfied. Let $u(t,x)$ be a classical solution of the Cauchy problem for the equation (1.3) with $f=0$ and the zero initial conditions. If for some $\sigma >0$,
\begin{equation}
\lim\limits_{|x|\to \infty}u(t,x)\exp \left\{ -\sigma |x|^{\frac2{2-\alpha}}\right\} =0
\end{equation}
uniformly with respect to $t\in [0,T]$, then $u(t,x)\equiv 0$.
\end{teo}

\medskip
Note that the formal substitution $\alpha =1$ corresponds to the classical uniqueness theorem for parabolic equations. The formal substitution $\alpha =2$ corresponds to the fact that no condition at infinity is needed for uniqueness for a hyperbolic equation. On the other hand, it is shown in \cite{Ps} for the equation (1.1) that the order $\dfrac2{2-\alpha}$ in (5.8) cannot be improved.

Note also that the above approach works also for the case $0<\alpha <1$ where the condition (5.8) and the assumption (D) guarantee the uniqueness too. For the equation (1.1) with $0<\alpha <1$ that is a result by Pskhu \cite{Ps}. For equations with variable coefficients, the uniqueness in a similar class of functions was known for $n=1$ (see \cite{K90,EK,EIK}); other uniqueness results for $0<\alpha <1$ \cite{EK,EIK,K90,K12} dealt with bounded solutions.

\section*{Acknowledgement}

This work was supported in part by Grant No. 03-01-12 of the National Academy of Sciences of Ukraine under the program of joint research projects with Siberian Branch, Russian Academy of Sciences.


\medskip
\begin{thebibliography}{999}
\bibitem{CGL}
Ph. Cl\'ement, G. Gripenberg, and S.-O. Londen, Regularity properties of solutions of fractional evolution equations, {Lect. Notes in Pure and Appl. Math.}, vol. 215 (2001), 235--246.
\bibitem{CLS}
Ph. Cl\'ement, S.-O. Londen, and G. Simonett, Quasilinear evolutionary equations and continuous interpolation spaces, {\it J. Diff. Equat.} {\bf 196} (2004), 418--447.
\bibitem{E}
S. D. Eidelman, Parabolic Systems, North-Holland, Amsterdam, 1969.
\bibitem{EK}
S. D. Eidelman and A. N. Kochubei, Cauchy problem for fractional diffusion equations. {\it J. Diff. Equat.} {\bf 199} (2004), 211--255.
\bibitem{EIK}
S. D. Eidelman, S. D. Ivasyshen, and A. N. Kochubei. {\it Analytic Methods in the Theory of Differential and Pseudo-Differential Equations of Parabolic Type}, Birkh\"auser, Basel, 2004.
\bibitem{F}
A. Friedman, {\it Partial Differential Equations of Parabolic Type}, Prentice-Hall, Englewood Cliffs, NJ, 1964.
\bibitem{Fu}
Y. Fujita, Integrodifferential equation which interpolates the heat equation and the wave equation, {\it Osaka J. Math.} {\bf 27} (1990), 309--321.
\bibitem{Ha}
A. Hanyga, Multidimensional solutions of time-fractional diffusion-wave equations, {\it Proc. Roy. Soc. London, Ser. A} {\bf 458} (2002), 933--957.
\bibitem{KST}
A. A. Kilbas, H. M. Srivastava, and J. J. Trujillo, {\it Theory and Applications of Fractional Differential Equations}, Elsevier, Amsterdam, 2006.
\bibitem{K90}
A. N. Kochubei, Fractional-order diffusion, {\it Differential Equations} {\bf 26} (1990), 485--492.
\bibitem{K12}
A. N. Kochubei, Fractional-parabolic systems, {\it Potential Anal.} {\bf 37} (2012), 1--30.
\bibitem{K13}
A. N. Kochubei, Fractional-hyperbolic systems, {\it Fract. Calc. Appl. Anal.} {\bf 16} (2013), 860--873.
\bibitem{LP}
Li Kexue and Peng Jigen, Fractional abstract Cauchy problem, {\it Integral Equ. Oper. Theory} {\bf 70} (2011), 333--361.
\bibitem{LMP}
Yu. Luchko, F. Mainardi and Yu. Povstenko, Propagation speed of the maximum of the fundamental solution to the fractional diffusion-wave equation, {\it Comp. Math. Appl.} {\bf 66} (2013), 774--784.
\bibitem{M96}
F. Mainardi, The fundamental solution for the fractional diffusion-wave equation, {\it Appl. Math. Lett.} {\bf 9} (1996), 23--28.
\bibitem{MP}
F. Mainardi and G. Pagnini, The Wright functions as solutions of the time-fractional diffusion equation, {\it Appl. Math. Comput.} {\bf 141} (2003), 51--62.
\bibitem{M}
F. Mainardi, {\it Fractional Calculus and Waves in Linear Viscoelasticity}, Imperial College Press, London, 2010.
\bibitem{Mi}
C. Miranda, {\it Partial Differential Equations of Elliptic Type}, Springer, Berlin, 1970.
\bibitem{PP}
H. Petzeltov\'a and J. Pr\"uss, Global stability of a fractional partial differential equation, {\it J. Integral Equat. Appl.} {\bf 12} (2000), 323--347.
\bibitem{PC}
I. Podlubny and Y. Chen, Adjoint fractional differential expressions and operators, Proc. ASME 2007 Conf. IDET/CIE 2007, Las Vegas, Nevada, USA. Paper DETC 2007--35005, 6 pages.
\bibitem{Pr}
J. Pr\"uss, {\it Evolutionary Integral Equations and Applications}, Birkh\"auser, Basel, 1993.
\bibitem{Ps}
A. V. Pskhu, The fundamental solution of a diffusion-wave equation of fractional order, {\it Izvestiya: Math.} {\bf 73} (2009), 351--392.
\bibitem{SKM}
S. G. Samko, A. A. Kilbas, and O. I. Marichev, {\it Fractional Integrals and Derivatives: Theory and Applications},
Gordon and Breach, New York, 1993.
\bibitem{SW}
W. R. Schneider and W. Wyss, Fractional diffusion and wave
equations, {\it J. Math. Phys.} {\bf 30} (1989), 134--144.
\bibitem{VZ}
V. Vergara and R. Zacher, Lyapunov functions and convergence to steady state for differential equations of fractional order, {\it Math. Z.} {\bf 259} (2008), 287--309.
\end{thebibliography}
\end{document}